                    \def\version{12 February, 2024}                       %
\documentclass[reqno,11pt]{amsart}
 \usepackage{amsmath, amsthm, a4, latexsym, amssymb, mathtools, enumitem, xcolor, color, tikz-cd,todonotes}
\usepackage[unicode]{hyperref}
\usepackage{srcltx}

\def\@rmrk#1#2{\refstepcounter
    {#1}\@ifnextchar[{\@yrmrk{#1}{#2}}{\@xrmrk{#1}{#2}}}

\makeatletter\@addtoreset{equation}{section}\makeatother

 \newfont{\bfit}{cmbxti10 scaled 1200}

\renewcommand{\d}{{\rm d}}

 \newcommand{\eps}{\varepsilon}
 \newcommand{\supp}{{\rm supp}}
 \newcommand{\Prop}{{\rm Prop}}
 \newcommand{\Tube}{{\rm Tube}}

 \newcommand{\R}{\mathbb{R}}
\newcommand{\N}{\mathbb{N}}
 \newcommand{\Z}{\mathbb{Z}}
 \newcommand{\C}{\mathbb{C}}

 \newcommand{\E}{\mathbb{E}}
 \renewcommand{\P}{\mathbb{P}}

 \def\1{{\mathchoice {1\mskip-4mu\mathrm l} 
{1\mskip-4mu\mathrm l}
{1\mskip-4.5mu\mathrm l} {1\mskip-5mu\mathrm l}}}

 \newcommand{\tmu}{\widetilde{\mu}}
 
 \newcommand{\talpha}{\widetilde{\alpha}}
 \newcommand{\tbeta}{\widetilde{\beta}}

 \newcommand{\tf}{\widetilde{f}}
 
 \newcommand{\tMcal}{\widetilde{\mathcal M}}
 \newcommand{\xtilde}{\widetilde{\mathcal X}}

 \newcommand{\skrif}{{\mathcal F}}

 \newcommand{\Bcal}{{\mathcal B}}
 \newcommand{\Mcal}{{\mathcal M}}

 \newcommand{\skris}{{\mathcal S}}

\newcommand{\bD}{{\mathbf D}}

\newcommand{\bS}{{\mathbf S}}
\newcommand{\bX}{{\mathbf X}}
\newcommand{\bK}{{\mathbf K}}

\newcommand{\bx}{{\mathbf x}}

\renewcommand{\subsection}{\secdef \subsct\sbsect}
\newcommand{\subsct}[2][default]{\refstepcounter{subsection}
\vspace{0.15cm}
{\flushleft\bf \arabic{section}.\arabic{subsection}~\bf #1  }
\nopagebreak\nopagebreak}
\newcommand{\sbsect}[1]{\vspace{0.1cm}\noindent
{\bf #1}\vspace{0.1cm}}

\newenvironment{example}{\refstepcounter{theorem}
{\bf Example \thetheorem\ }\nopagebreak  }%
{\nopagebreak {\hfill\rule{2mm}{2mm}}\\ }

\newtheorem{theorem}{Theorem}[section]
\newtheorem{lemma}[theorem]{Lemma}
\newtheorem{cor}[theorem]{Corollary}
\newtheorem{prop}[theorem]{Proposition}

\newtheorem{definition}[theorem]{Definition}

\newcommand{\eproof}{\hfill \qed \vspace*{5mm}}

\newtheoremstyle{thm}{1.5ex}{1.5ex}{\itshape\rmfamily}{}
{\bfseries\rmfamily}{}{2ex}{}

\newtheoremstyle{rem}{1.3ex}{1.3ex}{\rmfamily}{}
{\itshape\rmfamily}{}{1.5ex}{}
\theoremstyle{rem}
\newtheorem{remark}{{\slshape\sffamily Remark}}[]

\refstepcounter{subsubsection}

\def\thebibliography#1{\section*{References}
  \list%
  {\arabic{enumi}.}
    {\settowidth\labelwidth{[#1]}\leftmargin\labelwidth
    \advance\leftmargin\labelsep
    \parsep0pt\itemsep0pt
    \usecounter{enumi}}
    \def\newblock{\hskip .11em plus .33em minus .07em}
    \sloppy                   
    \sfcode`\.=1000\relax}

\begin{document}
\title[Schur Multipliers of $C^*$-algebras, amenability and group-invariant percolation]
{\large Schur Multipliers of $C^*$-algebras, group-invariant compactification and applications to amenability and percolation}
\author[Chiranjib Mukherjee and Konstantin Recke]{}
\maketitle
\thispagestyle{empty}
\vspace{-0.5cm}

\centerline{\sc Chiranjib Mukherjee\footnote{Universit\"at M\"unster, Einsteinstrasse 62, M\"unster 48149, Germany, {\tt chiranjib.mukherjee@uni-muenster.de}}
and Konstantin Recke\footnote{Universit\"at M\"unster, Einsteinstrasse 62, M\"unster 48149, Germany, {\tt konstantin.recke@uni-muenster.de}}}
\renewcommand{\thefootnote}{}

\vspace{-0.5cm}
\centerline{\textit{Universit\"at M\"unster}}
\vspace{0.2cm}

\begin{center}
\version
\end{center}

\begin{quote}{\small {\bf Abstract:} Let $\Gamma$ be a countable discrete group. Given any sequence $(f_n)_{n\geq 1}$ of  $\ell^p$-normalized functions ($p\in [1,2)$), 
consider the associated positive definite matrix coefficients $\langle f_n, \rho(\cdot) f_n\rangle$ of the right regular representation $\rho$. 
We construct an orthogonal decomposition of the corresponding {\it Schur multipliers} on the reduced group $C^*$-algebra
or the uniform Roe algebra of $\Gamma$. We identify this decomposition explicitly via the limit points of the orbits $(\widetilde f_n)_{n\geq 1}$ in the group-invariant compactification of the quotient space 
constructed by Varadhan and the first author in \cite{MV16}. We apply this result and use positive-definiteness to provide two (quite different) characterizations of amenability of $\Gamma$ -- one via a variational approach and the other
using group-invariant percolation on Cayley graphs constructed by Benjamini, Lyons, Peres and Schramm \cite{BLPS99}. These results underline, from a new point of view to the best of our knowledge, the manner in which Schur multipliers capture geometric properties 
of the underlying group $\Gamma$. 
}
\end{quote}


\section{Introduction and main results}\label{intro}

\subsection{Background.}\label{sec-background}
Let $\Gamma$ be a countable discrete group and denote by 
\begin{equation}\label{def-PGamma}
P(\Gamma):= \bigg\{\varphi \colon \Gamma \to \C \ \colon \ \forall \mbox{ finite }F\subset \Gamma, \, \, [\varphi(s^{-1} t)]_{s,t\in F} \,\,\,\,\mbox{is a positive definite matrix}\bigg\}
\end{equation}
the set of {\it positive definite functions} on $\Gamma$. These functions are of fundamental importance in diverse branches of mathematics for various reasons. First, in 
geometric group theory, properties like {\em amenability}, {\em Kazhdan's property (T)} and the {\em Haagerup property} admit characterizations in terms of pointwise limit points of positive definite functions. Moreover, they provide a well-known link to operator algebras via the following observation due to Haagerup \cite{Haagerup78}: {\it Schur multiplication} by the infinite matrix 
$$
[\varphi(st^{-1})]_{s,t\in \Gamma}, \quad\mbox{for}\,\,\varphi \in P(\Gamma),
$$
defines a {\it completely positive} linear operator $M_\varphi$ on the $C^*$-algebra $B(\ell^2(\Gamma))$ of bounded linear operators on $\ell^2(\Gamma)$, which is called a 
{\it Schur multiplier} (see Sec. \ref{sec-Schur}-\ref{sec-PosSchur}). Examples of such maps in operator algebras include \cite{BrownOzawa,GK02,O00,Roe,Y00} (see also \cite{B77,DD05,Pisier1996} for further background on Schur multipliers). Broadly speaking, finding suitable positive definite functions and analyzing the limit points of the associated Schur multipliers in the strong operator topology, plays a key role in 
understanding certain approximation properties of $C^*$-algebras associated to $\Gamma$, such as the {\it reduced group $C^*$-algebra} $C^*_\lambda(\Gamma)$ or the {\it uniform Roe algebra} $C^*_u(\Gamma)$ (see Subsection \ref{sex-unifRoeAlgebra} for details).

Perhaps the most canonical way of constructing positive definite functions goes as follows: Consider the map
\begin{equation}\label{def-Phi}
\begin{aligned}
\Phi: {\ell}^2(\Gamma) \to P(\Gamma), \quad  f\mapsto \Phi(f)(s) &= \sum_{t\in \Gamma} f(t) \overline{f(ts)} \\
&=\langle f, \rho(s) f \rangle,
\end{aligned}
\end{equation} 
where $\rho: \Gamma \to \mathcal U(\ell^2(\Gamma))$ denotes the {\it right regular representation}. Then $\Phi$ associates to an $\ell^2$-function $f$ the positive definite function $\Phi(f)$. Note that the class of positive definite functions that can be expressed in this way spans the right-analogue of the {\it Fourier algebra} of $\Gamma$ as defined by Eymard \cite{E64}. This class is considerably rich, containing in particular all finitely supported positive definite functions by a well-known result of Godement. In light of the above discussion, several questions then arise quite naturally: 

\begin{itemize}
\item[1.]  Given a sequence of Schur multipliers of the form $M_{\varphi_n}$ with $\varphi_n=\Phi(f_n)$, can we expect limit points in the strong operator topology? 
Is there a natural way to identify the limit points?$^{\rm a}$\footnote{$^{\rm a}$The same questions apply also to the possible pointwise limit points of the positive definite functions $\varphi_n=\Phi_n(f)$ themselves.}

\medskip

\item[2.] In the same setup, can the limit points of the sequence $M_{\varphi_n}$ with $\varphi_n=\Phi(f_n)$ be described via or recover information about the input sequence $f_n$ itself? 

\medskip

\item[3.]  Does the description of the above convergence and of the limit points shed new light on approximation properties of the group such as amenability?

\medskip

\item[4.] For the specific example of amenability, are there interesting probabilistic examples of the $\ell^2$-functions which give us this approximation property? We consider this aspect to be particularly relevant since any such probabilistic construction might potentially provide ways to deal with related, and more elusive, properties like Kazhdan's property (T) and the Haagerup property.

\end{itemize}

To underline the relevance of these questions, we note that for a countable group $\Gamma$, the input sequence $(f_n)_n$ (say, bounded in $\ell^2$) need not have {\it any} limit points in $\ell^2(\Gamma)$ -- 
the mass of $f_n$ already could split into widely separated pieces due to the action of the ambient group $\Gamma$ (e.g.\ $f_n(s)=\frac 12 \delta_n(s) + \frac 12 \delta_{-n}(s)$, $s\in \Z$), 
or the mass could totally disintegrate into dust (e.g.\ $f_n$ could be uniformly distributed on $\{-n,\dots, n\}$ on $\Z$). For similar reasons, possible  
limit points of the sequences $\{\Phi(f_n)\}_n$ or $\{M_{\Phi(f_n)\}_n}$ may not be tractable.

With this background, the goal and the guiding principle of our article will be to develop a unified approach to answer the above questions in a satisfactory manner. 
In fact, our first main result, stated in Theorem \ref{theorem1} below, will provide affirmative answers to the first two questions as corollaries. This is accomplished by developing a robust relationship between 
{\it decompositions} as well as convergence of Schur multipliers of the above form and convergence of the corresponding input functions $f$ in a {\em different space} -- namely, we will use 
the inherent group-invariance and structural properties of Roe algebras to leverage, in a unified way, 
the {\it group-invariant compactification} of the quotient space of orbits under the left $\Gamma$-action, we refer to Section \ref{sec-results1} below and the discussion following Theorem \ref{theorem1} for details. We then apply this theorem to prove our second main result, Theorem \ref{theorem3}, which uses the resulting compactness of the Schur multipliers to give a characterization of amenability and provides quantitative information on how well these Schur multipliers approximate the identity. This result may be considered to affirmatively answer the third question. Regarding the fourth question, we obtain a large class of interesting examples by employing a probabilistic model known as {\em group-invariant percolation} -- here we assume that $\Gamma$ is finitely generated. Informally speaking, we consider random subgraphs of some fixed Cayley graph of $\Gamma$ and study $\Phi(f)$ for random functions $f$ canonically arising from such models. This leads to our third main result, Theorem \ref{theorem4}, which is a second and quite different characterization of amenability. The approach also seems to be new and of independent interest, as we will elaborate on in more detail in the discussion following Theorem \ref{theorem4}. Let us now turn to a precise mathematical layout 
of the results alluded to in the above discussion.

\subsection{Main results.}\label{sec-results1}

We start with the function $\Phi$ defined in \eqref{def-Phi} and consider the following subsets of $\Phi(\ell^2(\Gamma))$: For $p\in[1,2]$, we set
\begin{equation}\label{def-Pp-Mcal}
\begin{aligned}
&\Mcal_p^{\scriptscriptstyle{\leq 1}}(\Gamma):= \big\{f\colon \Gamma \to [0,1]\colon \|f\|_p \leq 1\big\}, \quad
\Mcal_p(\Gamma):= \big\{f\colon \Gamma \to [0,1]\colon \|f\|_p = 1\big\}, \,\,\mbox{and}\\
&\qquad\qquad\qquad P_p(\Gamma):= \Phi\big(\Mcal_p^{\scriptscriptstyle{\leq 1}}(\Gamma)\big). 
\end{aligned}
\end{equation}
We remark that the restriction to non-negative functions with norm bounded by $1$ is made only for notational simplicity -- our results will work similarly when $\Mcal_p^{\scriptscriptstyle{\leq 1}}(\Gamma)$ is replaced with any $\|\cdot\|_p$-bounded set. Note that the case $p=2$ is essentially the general setting introduced in \eqref{def-Phi}, while for decreasing $p$-values the setting becomes increasingly restrictive. Now an important observation is that
$$
\Phi(f)= \Phi(\delta_s * f) \qquad \mbox{for all} \, \, s \in \Gamma, 
$$
meaning that $\Phi(f)$ depends on $f$ {\it only} via its {\it orbit} $\widetilde f=\{\delta_s * f : s \in \Gamma\}$ under the left action of $\Gamma$. 

Due to this observation, it is quite natural to pass to the quotient space $\widetilde\Mcal_p(\Gamma)$ of $\Mcal_p(\Gamma)$ under the left action of $\Gamma$. This space can be compactified using the approach developed in \cite{MV16} for orbits of probability measures on $\R^d$. A variant of the latter result and an adaptation to our 
discrete setting will be shown in Theorem \ref{theorem-CompactificationMp} -- in particular, we obtain that for every $p\in[1,2]$, the space 
\begin{equation}\label{def-Xcalp}
\widetilde{\mathcal X}_p(\Gamma):= \bigg\{ \xi= \{\widetilde\alpha_i\}_{i\in I}\colon I \subset \N, \alpha_i \in \Mcal_p^{\scriptscriptstyle{\leq 1}}(\Gamma), \,\,\sum_{i\in I} \|\alpha_i\|_p^p \leq 1 \bigg\}
\end{equation}
defines a compactification of the quotient $\widetilde\Mcal_p(\Gamma)$. We note that, in the above definition, $I \subset \N$ signifies that the collection $\{\widetilde \alpha_i\}_{i\in I}$ runs overs an empty, finite or countable collection $I$. 

We also note that in the most general case $p=2$, the map
\begin{equation}\label{map-l2}
\widetilde \Mcal_2(\Gamma) \ni \tf \mapsto \Phi(f)\in P(\Gamma)
\end{equation}
is in general not continuous -- indeed, for the normalized indicator functions 
$$
f_n = \frac{1}{\sqrt{2n+1}} \1_{\{-n,\ldots,n\}} \in \Mcal_2(\Z),
$$
 we have that $\widetilde f_n$ converges to $\widetilde 0$ because its $\ell^2$-mass restricted to {\it any} $r$-ball, i.e. $\sup_{x\in\Z} \Vert f_{|B(x,r)} \Vert_2^2$, goes to zero as $n\to\infty$ for every fixed $r>0$, but $\Phi(f_n)$ converges pointwise to $1 \neq \Phi(0)$. This is a general phenomenon, namely the same is true for the $\ell^2$-normalized indicator functions of any F{\o}lner sequence of an amenable group $\Gamma$. The situation for $p\in [1,2)$ is {\em completely different}: In this case, convergence and limit points of positive definite functions of the form $\Phi(f)$ and of their associated Schur multipliers are strongly determined by convergence of the orbits of the functions $f$ in the compactification (\ref{def-Xcalp}).

To state our result precisely, we introduce some terminology: For a Banach space $X$, let $B(X)$ denote the space of bounded linear operators from $X$ to itself, which is equipped with the strong operator topology. We define an {\it asymptotically orthogonal decomposition} (see Definition \ref{def-orthogonal}) of a sequence of operators $(T_n)_n$ in $B(X)$ w.r.t.\ some given subset $\mathcal S \subset B(X)$, to be a norm-convergent sum $\sum_{i\in I} S_i$, where $S_i \in \mathcal S$, such that, along a subsequence, $T_n \to \sum_i S_i$ in the strong operator topology and 
$\|T_n\| \to \sum_i \|S\|$. Note in particular that being an asymptotically orthogonal decomposition in this sense is a much stronger property than being a limit point. With this definition, we are now in a position to formally state our first main result (see Theorem \ref{theorem-MainTheorem}).

\begin{theorem}[Decomposing Schur multipliers and group-invariant compactification]\label{theorem1}
Let $\Gamma$ be a countable discrete group and fix $p\in[1,2)$. Let $(f_n)_{n=1}^\infty$ be any sequence in $\Mcal_p(\Gamma)$ and consider the sequence $(M_n)_{n=1}^\infty$ of Schur multipliers $M_n \coloneqq M_{\Phi(f_n)}$ on $C_u^*(\Gamma)$
or on $C_\lambda^*(\Gamma)$. Then the set of asymptotically orthogonal decompositions of $(M_n)_n$ with respect to 
$$
\left\{ M_{\Phi(f)} : f \in \Mcal_{p}^{\scriptscriptstyle{\leq 1}}(\Gamma) \right\}
$$
is given by 
$$
\left\{ M_\xi \coloneqq \sum_{\talpha \in \xi } M_{\Phi(\alpha)} \colon \xi \text{ is a limit point of } (\tf_n)_n \text{ in } \xtilde_p(\Gamma) \right\}
\vspace{2mm}
$$ 
and contains all limit points of the sequence $(M_n)_n$ in the strong operator topology (the above sum over the empty collection $\xi\in \xtilde_p$ is understood as zero).  Moreover, whenever the limit 
$$
\lim_{n \to \infty} M_{n} = M
$$
exists in the strong operator topology, every limit point $\xi \in \xtilde_p(\Gamma)$ of $(\tf_{n})_n$ must satisfy 
$$
M= M_\xi= \sum_{\widetilde\alpha\in \xi} M_{\Phi(\alpha)}.
$$
\end{theorem}

As mentioned previously, Theorem \ref{theorem1} actually entails several results as well as a number of consequences. To elaborate on that point, we observe that it directly yields affirmative answers to the first two questions raised in Section \ref{sec-background}. Indeed, it establishes a natural way of identifying the limit points of a sequence of Schur multipliers of the form $M_{\Phi(f_n)}$: namely, every limit point is given by a norm-convergent sum of Schur multpiliers of the same form, where the sum naturally corresponds to a limit point of the sequence $\tilde f_n$ in the compactification $\xtilde_p$ introduced above. Moreover, whether the sequence of Schur multipliers converges and to which limit, may be read off from knowing whether $\tilde f_n$ converges and to which limit. In other words, we have successfully transferred the problem of convergence of these Schur multipliers to the task of understanding convergence of the input functions in the compactification of the quotient space. One of the main advantages of this description is that convergence in the space $\xtilde_p$ actually has a very intuitive description: namely, a sequence $\tilde f_n$ converges to an empty, finite or countable collection $\xi=\{\tilde \alpha_i\}_i$ of orbits 
if the mass of $f_n$ comes from mutually widely separated pieces (in the metric of the underlying group $\Gamma$) $\alpha_i$, with some mass possibly dissipating to zero 
(see Section \ref{subsection: compactificationMcalp} for more details). Therefore, we also obtain a natural description of convergence of the associated Schur multipliers, which is formalized in our notion of asymptotically orthogonal decompositions. As an immediate and crucial consequence, we obtain compactness results for the sets of Schur multipliers introduced above. These clearly allow us to find maximizers and minimizers of suitable functionals and in particular of functionals appearing in the study of amenability. This leads to our second main result and will be explained in more detail in the next subsection. Let us also point out that similar results hold for positive definite functions of the form $\Phi(f)$ and seem to be new in this setting as well, see also Corollary \ref{corollary-limitpoints}. For the sake of completeness, we state the aforementioned consequences as corollaries below.

\begin{cor} Let $\Gamma$ be a countable discrete group and $p\in [1,2)$. Then 
$$
\bS_p(\Gamma) = \left\{ M_{\Phi(f)} \colon f \in \Mcal_p(\Gamma) \right\} \subset B(C^*_u(\Gamma))
$$ 
is relatively sequentially compact with respect to the strong operator topology and its closure in the strong operator topology is given by 
$$
\bX_p(\Gamma)=\{ M_\xi \colon \xi \in \xtilde_p(\Gamma) \}.
$$
\end{cor}

In the context of approximation properties, we are typically interested in approximating the identity by Schur multipliers or the constant function $1$ by positive definite functions. In the set-up considered above, we obtain the following consequence of compactness.

\begin{cor}\label{cor2}  Let $\Gamma$ be a countable discrete group and $p\in[1,2)$. Then for any finite subset $F\subset\Gamma$, there exist $f_i \in \Mcal_p^{\scriptscriptstyle{\leq 1}}, i \in \N$, with $\sum_{i=1}^\infty \|f_i\|_p^p \leq1$ such that the positive definite function
$$
\varphi=\sum_{i=1}^\infty \Phi(f_i)
$$
minimizes the expression
$$
d_F(\psi,1) \coloneqq \sup \bigl\{ |\psi(s)-1| \colon s \in F \bigr\}
$$
over all positive definite functions $\psi$ of this form and in particular over all $\psi \in P_p(\Gamma)$.
\end{cor}
On a related note, while strong operator topology convergence is the desired mode of convergence in the context of approximation properties of $C^*$-algebras, it is certainly natural to ask in which cases the mode of convergence in Theorem \ref{theorem1} can be improved to norm convergence. For this purpose, our notion of asymptotically orthogonal decompositions turns out to be a useful tool, using which we obtain the following, perhaps surprising, characterization of norm convergence (see Theorem \ref{thm-normconvergence}).

\begin{theorem}\label{theorem2} Let $\Gamma$ be a countable discrete group and $p\in[1,2)$. For any sequence $(f_n)_{n=1}^\infty$ in $\Mcal_p(\Gamma)$, consider the sequence $(M_n)_{n=1}^\infty$ of Schur multipliers $M_n\coloneqq M_{\Phi(f_n)}$ on $B(\ell^2(\Gamma))$. Then $(M_n)_n$ converges in norm to $M_{\Phi(\alpha)}$ if and only if there exists $\alpha\in\Mcal_p^{\scriptscriptstyle{\leq 1}}(\Gamma)$ such that
$$
\lim_{n\to\infty} \widetilde f_n = \{\widetilde \alpha\} \quad \mbox{in} \ \ \xtilde_p(\Gamma).
$$
\end{theorem}

The above statement means that $M_n$ converges in norm if and only if $\widetilde f_n$ converges to an element $\xi=\{\widetilde\alpha\}$ consisting of a {\it single orbit} $\widetilde\alpha$ in $\xtilde_p(\Gamma)$.

\subsection{Application to amenability and a connection with percolation theory.} 

We turn to our application of Theorem \ref{theorem1} to amenability. While the map in \eqref{map-l2} is not continuous for the case $p=2$, as remarked earlier, 
Theorem \ref{theorem1} applies for any $p \in [1,2)$ and we will show that it provides a robust approximation to reach the 
``threshold" $p=2$. To state this result, we introduce the following terminology: For a finite subset $F \subset \Gamma$, define $\Tube(F) \coloneqq \{ (x,y) \in \Gamma \times \Gamma \, \colon \, xy^{-1} \in F \}$ and let $M$ be a Schur multiplier. By setting all matrix entries outside $\Tube(F)$ to zero, we can make sense of the multiplier {\it cut-off outside $\Tube(F)$}, which will be denoted $r_F(M)$. 
Here is our next main result (contained in Theorem \ref{theorem: CompactApproximationId} stated later).

\begin{theorem}[A variational characterization of amenability using Schur multipliers]\label{theorem3}
Let $\Gamma$ be a countable discrete group. Then $\Gamma$ is amenable if and only if one of the following equivalent conditions hold:
\begin{enumerate}
\item The identity lies in the closure of 
$$
\bigcup_{p\in[1,2)} \bS_p(\Gamma)
$$
w.r.t.\ to the strong operator topology on $B(C_\lambda^*(\Gamma))$, respectively on $B(C_u^*(\Gamma))$. 
\item For each finite subset $F$ of $\Gamma$
\begin{equation} \label{eq-deltaFp}
\delta^{(F,p)} \coloneqq \min_{M \in \bX_p(\Gamma)} \Vert M_{\1_{\Tube(F)}} - r_F(M) \Vert 
\end{equation}
converges to zero as $p \in [1,2)$ converges to $2$. In this case, for any sequence of finite sets $F_n \uparrow \Gamma$, there exist $p_n \in [1,2)$ such that the minimizers $M^{(F_n,p_n)}$ in (\ref{eq-deltaFp}) converge, as $n\to\infty$, to the identity w.r.t.\ to the strong operator topology on  $B(C_\lambda^*(\Gamma))$, respectively on $B(C_u^*(\Gamma))$. 
\end{enumerate}
\end{theorem}
\noindent Theorem \ref{theorem3} asserts that for increasing finite subsets $F_n$ of our group $\Gamma$ and exponents $p_n$ increasing to $2$, it is possible to choose for each $n$ a {\em minimizer} in \eqref{eq-deltaFp} to obtain a sequence of Schur mulitpliers approximating the identity in the strong operator topology. Similarly, the corresponding positive definite functions converge pointwise to $1$. We stress that for fixed $p\in[1,2)$,  actually there is an optimal choice, which for $p=2$ typically can not exist. Roughly speaking, Theorem \ref{theorem3} characterizes amenability via solutions of minimization problems in the spaces introduced above, which were shown to be compact in Theorem \ref{theorem1}. In this sense, it incorporates much more information than just knowing whether the identity can be approximated by suitable maps and answers the third question raised in Section \ref{sec-background}.

\medskip

\noindent In the final part of this paper, we shift our attention to the setting of a finitely generated group $\Gamma$. We explore a new way of finding suitable sequences of positive definite functions using the probabilistic model of {\em percolations}, which are random subgraphs of some fixed Cayley graph of $\Gamma$ (see Section \ref{sec-percolation} for details).  More precisely, we investigate the behavior of positive definite functions of the form $\Phi(f)$, where $f$ is the normalized indicator function of an appropriate {\em finite} random subgraph. This leads to our final main result (see Theorem \ref{theorem: PercolationConstruction}).

\begin{theorem}[Amenability, percolations and positive definite functions]\label{theorem4}
Let $\Gamma$ be a finitely generated group, let $X$ be one of its Cayley graphs w.r.t.\ some finite and symmetric generating set and let $(p_n)_{n=1}^\infty \subset [1,2)$ such that $p_n \to 2$. Then $\Gamma$ is amenable if and only if there exists a sequence $(\P_n)_{n=1}^\infty$ of $\Gamma$-invariant site percolations on $X$ such that each $\P_n$ has no infinite components and for every $s\in\Gamma$
$$
\lim_{n\to\infty} \E_n\left[ \Phi \left( |K(o)|^{-1/p_n} \1_{K(o)} \right)(s) \right] = 1,
$$
where $K(o)$ denotes the random cluster of the identity in $X$ and $|K(o)|^{-1/p_n} \1_{K(o)} \coloneqq 0$ if $K(o)=\emptyset$.
Moreover, in this case 
$$
\varphi_n(s)\coloneqq\E_n\left[ \Phi \left( |K(o)|^{-1/p_n} \1_{K(o)} \right)(s) \right]
$$
 is positive definite and the Schur multipliers $M_{\varphi_n}$ on $B(C_\lambda^*(\Gamma))$, respectively on $B(C_u^*(\Gamma))$, converge to the identity in the strong operator topology as $n\to\infty$.
\end{theorem}

The above result can be shown using the {\it group-invariant percolations} constructed in \cite{BLPS99}  and the {\it mass-transport principle} introduced in \cite{H97} and further developed in \cite{BLPS99}. Its relevance in our context 
can be attributed to the fact that it allows to construct many examples of positive definite functions and Schur multipliers probabilistically using invariant percolations. 

Additional justification for the relevance of Theorem \ref{theorem4} is provided by the fact that it fits into a much broader context, which we now describe: The proof of the amenability criterion in \cite{BLPS99}, which inspired our result, uses the idea of {\em approximating an invariant mean} by expected averages over finite percolation clusters. Now for the purpose of understanding related properties like the Haagerup property or Kazhdan's property (T) using invariant percolations, this strategy (i.e., using approximate invariant mean) does not have an analogue. In contrast, we show that the same finite percolation clusters, informally speaking, also capture the information about positive definite functions which are relevant to amenability -- in this sense the class of {\em positive definite functions constructible via invariant percolations} is large enough to capture interesting properties of the group. This idea seems to have been overlooked in the literature so far. In fact, our approach of constructing positive definite functions via invariant percolations in Theorem \ref{theorem4} has a natural analogue in the context of the Haagerup property and property (T) -- this served as the main motivation for the percolation theoretic characterizations of these two properties which we recently obtained in \cite{MR23}.

\noindent{\bf Organization of the rest of the article:} In Section \ref{background} we will introduce Schur multipliers on operator algebras and provide the necessary facts; and in Section \ref{compactification} we will construct the aforementioned group-invariant compactification $\widetilde{\mathcal X}_p$ of $\ell^p$-bounded functions. Section \ref{main} constitutes our first main result Theorem \ref{theorem1}, its proof and the previously 
discussed consequences. Finally, Section \ref{section-applications} contains the aforementioned amenability criterion using Schur multipliers and the amenability criterion involving positive definite functions constructed from percolations. In particular, we will show Theorem \ref{theorem3} and Theorem \ref{theorem4} in Section \ref{sec-application1} and 
Section \ref{sec-percolation}, respectively.

\section{Schur multipliers of operator algebras}\label{background}

Throughout this section let $\Gamma$ be a countable discrete group. Consider the usual Hilbert space $\ell^2(\Gamma)$ of square-summable functions on $\Gamma$ with canonical orthonormal basis $\{ \delta_x \colon x \in \Gamma \}$ of Dirac masses and let $B(\ell^2(\Gamma))$ denote the $C^*$-algebra of bounded linear operators on $\ell^2(\Gamma)$. The following definitions are standard, see for example \cite{BrownOzawa, CW04, RW14}.

\subsection{Matrix representations.} \label{sec-matrixrepresentations}

Every operator $T \in B(\ell^2(X))$ can be uniquely represented as a $\Gamma \times \Gamma$ matrix $[T(x,y)]_{(x,y) \in \Gamma \times \Gamma}$ with coefficients defined by
$$
T(x,y) \coloneqq \langle \delta_x, T\delta_y \rangle.
$$
In this way $T$ can be viewed as an element of $\ell^\infty(\Gamma \times \Gamma)$ and we will write $\Vert T \Vert_\infty$ for the corresponding sup-norm. By the Cauchy-Schwarz inequality $\Vert T \Vert_\infty \leq \Vert T \Vert$, where $\Vert T \Vert$ is the operator norm of $T$. We say that $T$ has {\em constant diagonals} if $T(xs,ys)=T(x,y)$ for all choices $s,x,y \in \Gamma$. The set
$$
\supp (T) \coloneqq \{ (x,y) \in \Gamma \times \Gamma \, \colon \, T(x,y) \neq 0 \} \subset \Gamma\times\Gamma
$$
is called the {\em support} of $T$. For a subset $F \subset \Gamma$ define
$$
\Tube(F) \coloneqq \{ (x,y) \in \Gamma \times \Gamma \, \colon \, xy^{-1} \in F \}.
$$ 
The operator $T$ has {\em finite propagation} if there exists a finite subset $F \subset \Gamma$ such that
$$
\supp(T) \subset \Tube(F).
$$
In this case, the number
$$
\Prop(T) \coloneqq \inf \bigl\{ |F| \colon F \subset \Gamma, \supp(T) \subset \Tube(F) \bigr\}
$$
is called the {\em propagation} of $T$.

Conversely, given some matrix $[T(x,y)]_{(x,y)\in\Gamma\times\Gamma}$ there is no simple necessary and sufficient condition to decide whether a bounded linear operator on $\ell^2(\Gamma)$ with this matrix representation exists. We will use the following two well-known sufficient conditions.

\begin{lemma}[{\cite[Lemma 8.1]{WZ18}}] \label{lemma: FinPropBound} 
Let $[T(x,y)]_{(x,y)\in \Gamma \times \Gamma}$ be a matrix. Suppose there exists a finite subset $F$ of $\Gamma$ such that
$$
\bigl\{ (x,y) \in \Gamma \times \Gamma \colon T(x,y) \neq 0 \bigr\} \subset \Tube(F)
$$ 
and suppose that
$$
\sup_{(x,y) \in \Gamma \times \Gamma} |T(x,y)| = M < \infty.
$$ 
Then there exists a unique bounded linear operator $T$ on $\ell^2(\Gamma)$ with matrix representation $[T(x,y)]_{(x,y)}$. Moreover $\Vert T \Vert \leq |F| M.$ 
\end{lemma}

\proof We include a proof for the readers convenience. Let $\eta \in \ell^2(\Gamma)$. We observe that $T(x,y) \neq 0$ implies $xy^{-1} \in F$ which is equivalent to $y \in F^{-1}x = \{ y^{-1}x \colon y \in F\}$. Therefore for every $x \in \Gamma$ the sum
$$
\gamma_x \coloneqq \sum_{y \in \Gamma} T(x,y) \eta(y) = \sum_{y \in F^{-1}x} T(x,y)\eta(y)
$$
exists as it contains only finitely many summands. Moreover, by the Cauchy Schwarz inequality
$$
|\gamma_x|^2 \leq \Bigl( \sum_{y \in F^{-1}x} |T(x,y)|^2 \Bigr) \Bigl( \sum_{y \in F^{-1}x} |\eta(y)|^2 \Bigr) \leq |F| M^2  \sum_{y \in F^{-1}x} |\eta(y)|^2,
$$
which implies that
$$
\sum_{x \in \Gamma} |\gamma_x|^2 \leq |F|M^2\sum_{x \in \Gamma} \sum_{y \in F^{-1}x} |\eta(y)|^2 = |F|^2 M^2\Vert \eta \Vert_2^2.
$$
Hence $\eta \mapsto T\eta$, where $T\eta$ is defined by $(T\eta)(x) \coloneqq \gamma_x$, defines a bounded linear operator on $\ell^2(\Gamma)$ with operator norm $\Vert T \Vert \leq |F| M$ as claimed.
\eproof

\begin{remark} Let $\C_u[\Gamma]$ denote the collection of all operators $T\in B(\ell^2(\Gamma))$ with finite propagation. By Lemma \ref{lemma: FinPropBound}, $\C_u[\Gamma]$ can be identified with the collection of matrices with uniformly bounded entries and support contained in the tube of a finite set. Moreover, every $T \in \C_u[\Gamma]$ satisfies $\Vert T \Vert \leq \Prop(T) \Vert T \Vert_\infty$.
\end{remark}

\begin{lemma}[The weighted Schur test] \label{lemma: SchurTest} Let $[T(x,y)]_{(x,y) \in \Gamma \times \Gamma}$ be a matrix. Suppose there exists a function $p \colon \Gamma \longrightarrow (0,\infty)$ and two constants $C_1$ and $C_2$ such that for all $x,y \in \Gamma$
\vspace{1mm}
$$
\sum_{z \in \Gamma} |T(x,z)|p(z) \leq C_1p(x) \quad \text{and} \quad \sum_{z \in \Gamma} |T(z,y)|p(z) \leq C_2p(y).
$$
Then there exists a unique bounded linear operator $T$ on $\ell^2(\Gamma)$ with matrix representation $[T(x,y)]_{(x,y)}$. Moreover 
$$
\Vert T \Vert \leq (C_1C_2)^{\frac{1}{2}}. 
$$
\end{lemma}

\proof  
There are several ways to prove this result, see e.g.\ \cite[Section 3]{DK18} or \cite{Pisier1996}. For the readers convenience we include a direct proof similar to the proof of Lemma \ref{lemma: FinPropBound} above. To that end, let $\eta \in \ell^2(\Gamma)$ with finite support. Then for every $x\in\Gamma$, the sum
$$
\gamma_x \coloneqq \sum_{y \in \Gamma} T(x,y) \eta(y)
$$
exists as it contains only finitely many summands. Moreover, applying the Cauchy-Schwarz inequality in the fourth step and rearranging the appearing finite sums in the penultimate step, we obtain
\begin{flalign*}
\sum_{x\in\Gamma} |\gamma_x|^2 & = \sum_{x\in\Gamma} \Bigl| \sum_{y\in\Gamma} T(x,y) \eta(y) \Bigr|^2 \leq  \sum_{x\in\Gamma} \Bigl( \sum_{y\in\Gamma} |T(x,y)| |\eta(y)| \Bigr)^2\\
& = \sum_{x\in\Gamma} \Bigl( \sum_{y\in\Gamma} \bigl(|T(x,y)|p(y)\bigr)^{\frac{1}{2}} \biggl( \frac{|T(x,y)|}{p(y)}\biggr)^{\frac{1}{2}}  |\eta(y)| \Bigr)^2 \\
& \leq \sum_{x\in\Gamma} \Bigl( \sum_{y\in\Gamma} |T(x,y)| p(y) \Bigr) \Bigl( \sum_{y\in\Gamma} \frac{|T(x,y)|}{p(y)} |\eta(y)|^2 \Bigr) \\
& \leq \sum_{x\in\Gamma} C_1 p(x) \sum_{y\in\Gamma} \frac{|T(x,y)|}{p(y)} |\eta(y)|^2 \\
& = C_1 \sum_{y\in\Gamma} |\eta(y)|^2 \frac{1}{p(y)} \sum_{x\in\Gamma} |T(x,y)|p(x) \\
& \leq C_1 C_2 \Vert \eta \Vert_2^2.
\end{flalign*}
Since finitely supported elements are dense in $\ell^2(\Gamma)$, the claim follows. \eproof

\subsection{The uniform Roe algebra.}\label{sex-unifRoeAlgebra}

We recall the definitions of two important $C^*$-algebras associated with $\Gamma$ which are naturally represented on $B(\ell^2(\Gamma))$. First observe that it is not difficult to see that $\mathbb C_u[\Gamma]$, the collection of all finite propagation operators, defines a $^*$-subalgebra of $B(\ell^2(\Gamma))$. Its norm closure is called the {\em uniform Roe algebra} and is denoted by $C^*_u(\Gamma)$. The {\em reduced group $C^*$-algebra} $C_\lambda^*(\Gamma)$ is the $C^*$-algebra generated by the image of the {\em left-regular representation}
 $$
\lambda \colon \Gamma \longrightarrow B(\ell^2(\Gamma)) \, , \, s \mapsto \lambda(s),
$$
where 
$$
(\lambda(s)f)(x) = f(s^{-1}x) \quad \text{for } f \in \ell^2(\Gamma).
$$ 
Equivalently the reduced group $C^*$-algebra is the norm closure of the $^*$-subalgebra of $B(\ell^2(\Gamma))$ which consists of all operators $T \in B(\ell^2(\Gamma))$ with finite propagation and constant diagonals. We also recall that the {\em right-regular representation} of $\Gamma$ is the map
$$
\rho \colon \Gamma \longrightarrow B(\ell^2(\Gamma)) \, , \, s \mapsto \rho(s),
$$ 
where
$$
(\rho(s)f)(x) = f(xs) \quad \text{for } f \in \ell^2(\Gamma).
$$

\begin{remark} The group $\Gamma$ has a natural left action on the set $\Mcal_1(\Gamma)$ of probability measures on $\Gamma$ via the maps $\mu \mapsto \mu(s^{-1} \cdot) = \delta_s * \mu$, see the beginning of Section \ref{compactification} for more details. This action is implemented by the left-regular representation in the following sense: If $\mu \in \Mcal_1(\Gamma)$, let $f \in \ell^1(\Gamma) \subset \ell^2(\Gamma)$ denote its unique density with respect to counting measure (Haar measure). Then for every $s,x\in\Gamma$
$$
(\lambda(s)f)(x)=f(s^{-1}x)= \sum_{y \in \Gamma} \delta_s(y)f(y^{-1}x) = (\delta_s * \mu)(\{x\}),
$$
that is $\delta_s * \mu$ has density $\lambda(s)f$.
\end{remark}

\subsection{Schur multipliers.}\label{sec-Schur}

For a function $k \in \ell^\infty(\Gamma \times \Gamma)$ and an operator $T \in B(\ell^2(\Gamma))$ the {\em Schur product} is defined as
$$
k \circ T \coloneqq [k(x,y)T(x,y)]_{(x,y) \in \Gamma \times \Gamma}.
$$
In general such a Schur product need not be the matrix representation of a bounded linear operator. In the case that $k \circ T$ is the matrix representation of a bounded linear operator for every $T \in B(\ell^2(\Gamma))$, we denote this operator also by $k \circ T$ and obtain a well-defined linear map 
$$
M_k \colon B(\ell^2(\Gamma)) \longrightarrow B(\ell^2(\Gamma)) \, , \, T \mapsto M_k(T) \coloneqq k \circ T,
$$
which is called a {\em Schur multiplier}. It then follows from the Closed Graph Theorem that every Schur multiplier is bounded.

\begin{example} Let $k \in \ell^\infty(\Gamma \times \Gamma)$ such that $\supp(k) \subset \Tube(F)$ for some finite subset $F$ of $\Gamma$. Then Lemma \ref{lemma: FinPropBound} shows that $M_k$ is a Schur multiplier. Moreover, we observe the following:
\begin{enumerate}[label=(\alph*)]
\item The map $M_k$ restricts to a well-defined bounded linear operator on $C_u^*(\Gamma)$.
\item If the function $k$ is {\em diagonally invariant}, in the sense that $k(xs,ys)=k(x,y)$ for all $s,x,y\in\Gamma$, then $M_k$ also restricts to a well-defined bounded linear operator on $C_\lambda^*(\Gamma)$.
\end{enumerate}
A typical example for (b) is when $k$ is of the form $k(x,y)=f(xy^{-1})$ for some function $f\in C_c(\Gamma)$.
\end{example}

We note that every function $f$ on $\Gamma$ gives rise to a natural diagonally invariant kernel
$$
k_f \colon \Gamma \times \Gamma \longrightarrow \mathbb C \, , \, k_f(x,y)=f(xy^{-1}).
$$
If the corresponding Schur multiplier exists, we will denote it by $M_f \coloneqq M_{k_f}$ to simplify notation. In this case $M_f$ restricts to a bounded linear operator on $C_u^*(\Gamma)$ as well as on $C_\lambda^*(\Gamma)$, which will be called the corresponding {\it Schur multiplier on} $C_u^*(\Gamma)$, resp.\  {\it on} $C_\lambda^*(\Gamma)$. If it is clear from the context which of these maps we are considering, the Schur mutliplier on $C_u^*(\Gamma)$, resp.\ $C_\lambda^*(\Gamma)$, will also be denoted $M_f$. The following basic lemma shows that these multipliers exist for example for all integrable functions.

\begin{lemma} \label{lemma: SchurMultBounded}  Let $f \in \ell^1(\Gamma)$. Then the Schur multiplier
$$
M_f \colon B(\ell^2(\Gamma)) \longrightarrow B(\ell^2(\Gamma)) \, , \, T \mapsto M_f(T) \coloneqq [k_{f}(x,y)T(x,y)]_{(x,y)\in \Gamma \times \Gamma},
$$
is a well-defined bounded linear operator with $\Vert M_f \Vert \leq \Vert f \Vert_1$. In particular the restriction of $M_f$ to $C_u^*(\Gamma)$, resp.\ $C_{\lambda}^*(\Gamma)$, is a bounded linear operator on $C_u^*(\Gamma)$, resp.\ $C_{\lambda}^*(\Gamma)$, with norm bounded by $\Vert f \Vert_1$.
\end{lemma}

\proof Let $T \in B(\ell^2(\Gamma))$ and set $S \coloneqq k_f \circ T = [ f(xy^{-1})T(x,y) ]_{(x,y) \in \Gamma \times \Gamma}$. For each $x\in\Gamma$
$$
\sum_{x \in \Gamma} |S(x,y)| = \sum_{x \in \Gamma} |f(xy^{-1})T(x,y)| = \sum_{z \in \Gamma} |f(z)T(zy,y)| \leq \Vert T \Vert_\infty \Vert f \Vert_1 \leq \Vert T \Vert \Vert f \Vert_1
$$
and similarly for each $y\in\Gamma$
$$
\sum_{y \in \Gamma} |S(x,y)| =  \sum_{y \in \Gamma} |f(xy^{-1})T(x,y)| = \sum_{z \in \Gamma} |f(z)T(x,z^{-1}x)| \leq \Vert T \Vert_\infty \Vert f \Vert_1 \leq \Vert T \Vert \Vert f \Vert_1.
$$
Applying Lemma \ref{lemma: SchurTest} with weights $p\equiv1$ we conclude that $S=M_f(T) \in B(\ell^2(\Gamma))$ with 
$$
\Vert S \Vert \leq \Vert T \Vert \Vert f \Vert_1 .
$$
It follows that $M_f$ defines a Schur multiplier and that $\Vert M_f \Vert \leq \Vert f \Vert_1$. 
\eproof

\subsection{Completely positive Schur multipliers.}\label{sec-PosSchur}

Finding the norm of a Schur multiplier is a difficult problem in general. However, it is possible in special cases \cite{DD05,HSS10}, for example if the Schur multiplier is induced by a positive definite function. To discuss this example, we briefly recall the necessary vocabulary: A function $\varphi \colon\Gamma\longrightarrow \mathbb C$ is called {\it positive definite} if the matrix $[\varphi(s^{-1}t)]_{s,t \in F} \in M_{|F|}(\mathbb C)$ is positive for every finite set $F \subset \Gamma$. For a $C^*$-algebra $A$ and $n\geq1$, let $\mathbb M_n(A)$ denote the $C^*$-algebra of $n\times n$ matrices with entries in $A$ equipped with the usual addition, matrix multiplication and involution. A function $\phi$ mapping a unital $C^*$-algebra $A$ to a $C^*$-algebra $B$ is {\it completely positive} if the induced map $\phi_n \colon \mathbb M_n(A) \longrightarrow \mathbb M_n(B), \phi_n([a_{i,j}])=[\phi(a_{i,j})]$ is positive for every $n\geq1$. It is {\it completely bounded} if the sequence $(\phi_n)_n$ is uniformly bounded and in this case
$$
\Vert \phi \Vert_{cb} \coloneqq \sup_{n\geq1} \Vert \phi_n \Vert
$$
is called the {\it cb-norm} of $\phi$.

\begin{lemma}[{\cite[Lemma 1.1]{Haagerup78}}] \label{lemma: Haagerup}
Let $\varphi$ be a positive definite function on $\Gamma$. Then 
$$
M_\varphi \colon C_\lambda^*(\Gamma) \longrightarrow C_\lambda^*(\Gamma)
$$
is completely positive with norm $\Vert M_\varphi \Vert = \varphi(e)$.
\end{lemma}

\proof
This is a well known result of Haagerup, for a proof see Lemma 1.1 in \cite{Haagerup78}. \eproof

Using standard properties of completely positive maps, it is not difficult to see that the above result holds as well for Schur multipliers on $C_u^*(\Gamma)$ and $B(\ell^2(\Gamma))$. For the convenience of the reader we will state this slight generalization as a lemma and include a proof.

\begin{lemma} \label{lemma: GeneralizedHaagerup} Let $\varphi$ be a positive definite function on $\Gamma$. Then
$$
M_\varphi \colon B(\ell^2(\Gamma)) \longrightarrow B(\ell^2(\Gamma))
$$
is completely positive with norm $\Vert M_\varphi \Vert = \varphi(e)$ and cb-norm $\Vert M_\varphi \Vert_{cb} =\varphi(e)$.  \\ 
Moreover, its restriction $M_\varphi^u$ to $C_u^*(\Gamma)$ and its restriction $M_\varphi^\lambda$ to $C_\lambda^*(\Gamma)$ are both completely positive with norm and cb-norm equal to $\varphi(e)$.
\end{lemma}

\proof 
It is a standard fact that every completely positive map is completely bounded and that its norm and cb-norm agree, see for example \cite[Appendix B, p.\ 450]{BrownOzawa}.  Thus Lemma \ref{lemma: Haagerup} already implies the statement about $M_\varphi^\lambda$.

By \cite[Proposition D.3]{BrownOzawa} the Schur multiplier $M_\varphi$ is completely positive and hence its norm and cb-norm coincide. Moreover, by \cite[Proposition D.6]{BrownOzawa}, $M_\varphi$ is completely bounded if and only if $M_\varphi^\lambda$ is completely bounded and in this case their cb-norms coincide. Hence the statement about $M_\varphi$ follows from the one about $M_\varphi^\lambda$.

Finally $\Vert M_\varphi^\lambda \Vert \leq \Vert M_\varphi^u \Vert \leq \Vert M_\varphi \Vert$ implies that $\Vert M_\varphi^u \Vert = \varphi(e)$. The fact that $M_\varphi^u$ is also completely positive follows for example from Stinespring's Dilation Theorem, see \cite[Theorem 1.5.3]{BrownOzawa}, and thus its norm and cb-norm coincide. Therefore we have proved the statement for $M_\varphi^u$ and the proof of Lemma \ref{lemma: GeneralizedHaagerup} is thus complete. \eproof

Let us also mention the following fundamental result which is far more general than Lemma \ref{lemma: SchurMultBounded}. It is originally due to Grothendieck \cite{Grothendieck}. The version below and a proof can be found in \cite{Pisier1996}.

\begin{theorem} \label{theorem: Grothendieck} 
Let $A$ and $B$ be arbitrary sets, $k\colon A \times B \longrightarrow \C$ be a function and $C\geq 0$ be a constant. Then the following are equivalent:
\begin{enumerate}
\item The map
$$ 
M_k \colon B(\ell^2(A),\ell^2(B)) \longrightarrow B(\ell^2(A),\ell^2(B))
$$
which maps a bounded linear operator $T\colon \ell^2(A) \longrightarrow \ell^2(B)$ with matrix $[T(x,y)]_{x\in A, y \in B}$, where again $T(x,y)=\langle \delta_x, T\delta_y \rangle$, to a bounded linear operator with matrix $[k(x,y)A(x,y)]$, is well-defined (and thus bounded) with norm $\Vert M_k \Vert \leq C$.
\item The map $M_k$ is completely bounded with cb-norm $\Vert M_k \Vert_{cb} \leq C$.
\item There exist a Hilbert space $H$ and two families $(\varphi_x)_{x \in A}$ and $(\psi_y)_{y \in B}$ in $H$ such that 
$$
\sup_{x\in A} \Vert \varphi_x \Vert \sup_{y\in B} \Vert \psi_y \Vert \leq C
$$ 
and
$$
k(x,y)=\langle \varphi_x, \psi_y \rangle \quad \text{for all } x \in A, y \in B.
$$
\end{enumerate}
\end{theorem}

\begin{remark} The above result implies for example that $M_f$ is a Schur multiplier with norm $\Vert M_f \Vert \leq \Vert f \Vert_2$ for every $f\in\ell^2(\Gamma)$. This follows by choosing the sets $A,B=\Gamma$, the Hilbert space $H=\ell^2(\Gamma)$ and the functions $\varphi_x(z)=f(xz^{-1})$ and $\psi_y=\delta_y$. However, we will mainly be concerned with Schur multipliers arising from positive definite functions for which the above theorem is not needed.
\end{remark}

\section{Group invariant compactification}\label{compactification}

We will start with the general setting and then focus on the case of countable discrete groups. Consider a Polish group $G$ with a compatible left-invariant metric $d_G$, i.e.\ $G$ is a topological group and $d_G$ is a left-invariant metric on $G$ which induces the topology and is such that $(G,d_G)$ is complete and separable. The Borel-$\sigma$-field on $G$ will be denoted by $\Bcal(G)$. We denote by $\Mcal_1 = \Mcal_1(G)$ the space of probability measures on $(G,\Bcal(G))$ and by $\Mcal_{\leq1} = \Mcal_{\leq 1}(G)$ the space of sub-probability measures, i.e.\ non-negative measures on $(G,\Bcal(G))$ with total mass less than or equal to one. We consider the natural left action of $G$ on $\Mcal_{\leq 1}$ via the maps 
$$
\alpha \mapsto \alpha(s^{-1} \cdot) = \delta_s * \alpha \quad \text{for } \alpha \in \Mcal_{\leq 1}, s \in G.
$$
For any $\alpha \in \Mcal_{\leq 1}$, its {\it orbit} under this action is denoted by $\talpha = \{ \delta_s*\alpha \colon s \in G \}$. We denote by $\tMcal_1$ the quotient space of $\Mcal_1$ and by $\tMcal_{\leq 1}$ the quotient space of $\Mcal_{\leq1}$ under this action.

\subsection{Topologies on measures.} 

We will now recall the main properties of the topologies which are needed to understand the compactification in the next subsection. More details can be found in \cite{MV16} and more background on convergence of probability measures may be found e.g.\ in \cite{Billingsley}. Let us first describe two natural topologies on $\Mcal_{\leq1}$. In the {\it weak topology} a sequence $\mu_n$ converges to $\mu$, denoted by $\mu_n \Rightarrow \mu$, if
\begin{equation} \label{equ: criterionweakconv}
\lim_{n\to\infty} \int_{G} f(x) \mu_n(dx) = \int_{G} f(x) \mu(dx)
\end{equation}
for every bounded continuous function $f$ on $G$. Contrary, in the {\it vague topology}, a sequence $\mu_n$ converges to $\mu$, denoted by $\mu_n \hookrightarrow \mu$, if the convergence in \eqref{equ: criterionweakconv} holds for all continuous functions with compact support. In this case it continues to hold for all continuous functions vanishing at infinity. If $G$ is compact, then weak and vague topology agree and $\Mcal_1$ equipped with this topology is compact. Therefore we will only be interested in the case of a non-compact group.

The weak topology is the usual topology from the point of view of probability theory, an important reason being that $\Mcal_1$ is closed in  $\Mcal_{\leq1}$ with respect to weak convergence and thus weak convergence preserves the total amount of mass. However, for non-compact groups the space $\Mcal_1$ is not compact in this topology. There are two typical reasons for this: The location of the mass may shift away to infinity, as in $\mu_n = \delta_{n} \in \Mcal_1(\R)$, or 
the mass may totally disintegrate into dust, as in the example of a Gaussian $\mu_n = N(0,\sigma_n) \in \Mcal_1(\R)$ with variance $\sigma_n$ diverging to infinity.
Of course one might also have a combination of these two situations, as in the example 
$\mu_n =  \frac{1}{3}\bigl(\delta_n + \delta_{n^2} + N(0,\sigma_n)\bigr)$,
where the mass splits into two pieces increasingly far away from one another and another piece which disintegrates into dust. The vague topology provides one possible solution to the lack of compactness in $\Mcal_1$ because $\Mcal_{\leq1}$ is a compactification with respect to the vague topology. The problem with this compactification is that it essentially amounts to ignoring the mass which disintegrated into dust as well as all mass which moved off to infinity. Such sequences will all converge vaguely to zero, although the two types of behaviors are fundamentally different.

\subsection{Compactification of $\tMcal_1$: The space $\xtilde$.} 

Given the above set up, the following idea was introduced in \cite{MV16} and a {\it compactification of the quotient space} $\widetilde\Mcal_1$ was constructed there. Namely, consider the quotient space $\tMcal_1(\R^d)$ and equip this space with the quotient topology. Because of the preceding discussion, $\tMcal_1(\R^d)$ need not be compact either. So a compactification of $\tMcal_1(\R^d)$ should then correspond to a space $\xtilde$ consisting of empty, finite or countable collections of orbits, whose total masses add up to at most one. In order to make this intuition precise, it is useful to define a totally bounded compatible metric on $\tMcal_1$ which can naturally be extended to $\xtilde$ and then prove that this space is a compactification. We will now describe this construction in a general set-up.

Define the following sets of test functions: For $k\geq 2$, let $\skrif_k$ be the set of all continuous functions
$$
f \colon G^k \longrightarrow \R,
$$
which are \textit{diagonally invariant} in the sense that
$$
f(sx_1,\ldots,sx_k)=f(x_1,\ldots,x_k) \quad \text{ for all } s,x_1,\ldots,x_k \in G,
$$
and \textit{vanish at infinity} in the sense that
$$
\lim_{\max_{i \neq j} d_G(x_i,x_j) \to \infty} f(x_1,\ldots,x_k) = 0.
$$
These two properties imply that $\skrif_k$ can be identified with the space $C_0\bigl(G^{k-1},\R\bigr)$ of continuous real valued functions of $(k-1)$ variables vanishing at infinity in the usual sense. Thus $\skrif_k$ equipped with the uniform metric is a separable metric space. In particular there exists a sequence $(f_r)_{r=1}^\infty \subset \bigsqcup_{k=2}^\infty \skrif_k$ containing for every $k\geq2$ a subsequence which is dense in the corresponding space $\skrif_k$. For each $r\in\N$ let $k_r\geq2$ be the corresponding index such that $f_r\in\skrif_{k_r}$.

Define
$$
\xtilde \coloneqq \xtilde(G) \coloneqq \Bigl\{ \xi = \{ \talpha_i \}_{i \in I} \, \colon \, I \subseteq \N, \talpha_i \in \widetilde \Mcal_{\leq1}, \sum_{i \in I} \alpha_i(G) \leq 1 \Bigr\}
$$
with the convention that the empty collection is identified with the set $\{ \tilde 0\}$ containing only the orbit of the zero measure. Then $\xtilde$ naturally contains the space $\widetilde \Mcal_1$ via the identification of $\tmu$ with $\{\tmu\}$. For $\xi,\xi'\in\xtilde$ define
\begin{equation}\label{def-D}
\begin{aligned}
\bD(\xi,\xi')\coloneqq \sum_{r=1}^\infty \frac{1}{2^r} \frac{1}{1+\Vert f_r\Vert_\infty} \Bigg| \sum_{\talpha \in \xi} \int & f_r(x_1,\ldots,x_{k_r}) \prod_{i=1}^{k_r} \alpha(\d x_i) \\ 
& - \sum_{\talpha \in \xi'} \int f_r(x_1,\ldots,x_{k_r}) \prod_{i=1}^{k_r} \alpha(\d x_i) \Bigg|
\end{aligned}
\end{equation} 
and note that diagonal invariance of the functions $f_r$ ensures that the above definition does not depend on the choice of orbit representatives. We are now in a position to recall the main result obtained in \cite{MV16} regarding the space $\xtilde$.

\begin{theorem} \label{theorem: MVcompactification} $(\xtilde(\R^d), \bD)$ is a compact metric space which contains $\widetilde \Mcal_1(\R^d)$ as a dense subspace. It is then also the completion under the metric $\bD$ of the totally bounded space $\tMcal_1(\R^d)$. 
\end{theorem}

\proof 
See Theorem 3.1 and Theorem 3.2 in \cite{MV16}. \eproof

\subsection{Countable discrete groups and $\ell^p$-bounded densities.} \label{subsection: compactificationMcalp}

For the purposes of this paper, we will now consider the setting of a countable discrete group $\Gamma$. In this context, let us take note of the following basic facts. 
First, note that every countable discrete group admits a proper left-invariant metric, which is also unique up to coarse equivalence, compare e.g.\  \cite[Proposition 5.5.2]{BrownOzawa}.  From this point on we will assume that $\Gamma$ is equipped with a fixed proper left-invariant metric $d_{\Gamma}$.
Next, every measure on $\Gamma$ has a density with respect to counting measure, i.e.\ Haar measure, and we may thus identify 
\begin{equation} \label{equ: M1discrete}
\Mcal_1(\Gamma) = \left\{ \mu \colon \Gamma \longrightarrow [0,1] \colon \Vert \mu \Vert_1 = 1 \right\} \subset \ell^1(\Gamma).
\end{equation}
Let $(\mu_n)_{n=1}^\infty$ and $\mu$ be in $\Mcal_{\leq1}(\Gamma)$. It is not difficult to show that $\mu_n \Rightarrow \mu$ if and only if $\Vert \mu_n - \mu \Vert_1 \to 0$ and that $\mu_n \hookrightarrow \mu$ if and only if $\mu_n \to \mu$ pointwise.
Analogous to the identification in \eqref{equ: M1discrete}, we define for each $p \in [1,\infty)$ the sets
$$
\Mcal_p \coloneqq \Mcal_p(\Gamma) \coloneqq \left\{ f \colon \Gamma \longrightarrow [0,1] \colon \Vert f \Vert_p = 1 \right\} \subset \ell^p(\Gamma),
$$
$$
\Mcal_p^{\scriptscriptstyle{\leq 1}}\coloneqq \Mcal_p^{\scriptscriptstyle{\leq 1}}(\Gamma) \coloneqq \left\{ f \colon \Gamma \longrightarrow [0,1] \colon \Vert f \Vert_p \leq 1 \right\} \subset \ell^p(\Gamma).
$$
We consider the same left action of $\Gamma$ and denote the corresponding quotient spaces by $\tMcal_p$ and $\tMcal_p^{\scriptscriptstyle{\leq 1}}$. 
Note that $\Mcal_p$ fails to be compact with respect to the norm topology and with respect to pointwise convergence for similar reasons as in the previous case of $\Mcal_1$. We will show in this subsection that the compactification $\xtilde$ induces a compactification of $\tMcal_p$, equipped with the quotient of the $\ell^p$-norm topology, which again recovers exactly the pieces of mass lost due to the action of the group. To that end, define for each $p\in(1,\infty)$
\vspace{1mm}  
$$
\xtilde_p = \xtilde_p(\Gamma) \coloneqq \Bigl\{ \{ \talpha_i \}_{i \in I} \, \colon \, I \subseteq \N, \talpha_i \in \tMcal_p^{\scriptscriptstyle{\leq 1}}, \sum_{i \in I} \Vert\alpha_i \Vert_p^p \leq 1 \Bigr\},
$$
with the convention that the empty collection is identified with the set $\{ \tilde 0\}$ containing only the orbit of the constant zero function. 
Now consider the map
$$
\begin{aligned}
&\Mcal_{\leq1} \longrightarrow \Mcal_p^{\scriptscriptstyle{\leq 1}} \\
& \alpha \mapsto \alpha^{1/p}
\end{aligned}
$$
and extend it in the natural way to $\xtilde$ by defining for $\xi=\{\talpha_i \}_{i\in I} \in \xtilde$
$$
\xi^{1/p} \coloneqq \bigl\{ \talpha_i^{1/p} \bigr\}_{i\in I} \coloneqq \biggl\{ \widetilde{\alpha_i^{1/p}} \biggr\}_{i\in I} \in \xtilde_p.
$$
To see that this map is well-defined, observe that
$$
\sum_{i \in I} \Vert \alpha_i^{1/p} \Vert_p^p = \sum_{i \in I} \Vert \alpha_i \Vert_1 = \sum_{i \in I} \alpha_i(\Gamma).
$$
In particular the map $\xi \mapsto \xi^{1/p}$ is onto. Since it is clearly injective, it is one-to-one and its inverse is given by 
$$
\xtilde_p \longrightarrow \xtilde \, , \, \{\tbeta_i \}_{i\in I}=\eta \mapsto \eta^p \coloneqq \bigl\{ \tbeta_i^p \bigr\}_{i \in I}.
$$
Recall the definition of $\bD$ from \eqref{def-D}. Then for any $\eta,\eta'\in\xtilde_p$ we now set
$$
\bD_p(\eta,\eta') \coloneqq \bD\bigl(\eta^p,(\eta')^p\bigr).
$$
By an argument similar to that of Theorem \ref{theorem: MVcompactification} it can be shown that $\bD_p$ is a metric. Finally, we verify that $\bD_p$ induces on $\tMcal_p$ the quotient of the norm-topology: Let $(f_n)_{n=1}^\infty$ and $f$ be in $\Mcal_p$ and define 
$$
\mu_n \coloneqq f_n^p\in\Mcal_1\quad\mbox{ and }\quad\mu \coloneqq f^p\in\Mcal_1.
$$
 By definition
$$
\bD_p (\tf_n,\tf) \to 0 \quad \text{if and only if} \quad \bD(\tmu_n,\tmu) \to 0.
$$
The latter conditions is equivalent to the existence of shifts $s_n$ such that $\delta_{s_n} * \mu_n \Rightarrow \mu$. Now on the one hand, if 
$$
\Vert \delta_{s_n}*f_n-f\Vert_p \to 0,
$$
 then clearly $\delta_{s_n}*\mu_n$ converges pointwise to $\mu$ and since $\delta_{s_n}*\mu_n(\Gamma)\equiv1$ and $\mu(\Gamma)=1$, a standard argument using tightness implies that $\delta_{s_n} * \mu_n \Rightarrow \mu$. Conversely, if $\delta_{s_n} * \mu_n \Rightarrow \mu$, then clearly $\delta_{s_n} * f_n$ converges pointwise to $f$ and since $\Vert \delta_{s_n}*f_n \Vert_p \equiv 1$ and $\Vert f \Vert_p = 1$, a similar argument shows that $\Vert \delta_{s_n}*f_n-f\Vert_p \to 0$. 
 We will now show the following analogue of Theorem \ref{theorem: MVcompactification}.

\begin{theorem}\label{theorem-CompactificationMp} Let $\Gamma$ be a countable discrete group and $p\in[1,\infty)$. Then \label{CompactificationMp} $(\xtilde_p(\Gamma), \bD_p)$ is a compact metric space which contains $\widetilde \Mcal_p(\Gamma)$ as a dense subspace. It is then also the completion under the metric $\bD_p$ of the totally bounded space $\tMcal_p(\Gamma)$.
\end{theorem}

\subsection{\bf Proof of Theorem \ref{theorem-CompactificationMp}.}

The proof follows the line of argument developed in \cite[Proof of Theorem 3.2]{MV16} with an appropriate analogue of the concentration function of a probability measure. More specifically, we define the {\it concentration function} $q_f$ of $f\in\Mcal_p^{\scriptscriptstyle{\leq 1}}$ by
$$
q_f(r) \coloneqq \sup_{x\in\Gamma} \big\Vert f |_{B(x,r)} \big\Vert_p^p \qquad (r>0).
$$ 
Note that
$$
\lim_{r\to\infty} q_f(r)=\Vert f \Vert_p^p.
$$
The decomposition now arises from the following iterative procedure: Start with any sequence $(f_n)_{n=1}^\infty$ in $\Mcal_p^{\scriptscriptstyle{\leq 1}}$ and write $q_n\coloneqq q_{f_n}$ to shorten notation. By passing to a subsequence if necessary, we may assume that 
$$
q(r) \coloneqq \lim_{n\to\infty} q_n(r)
$$
exists for every $r>0$ and that 
$$
p^* \coloneqq \lim_{n\to\infty} \Vert f_n \Vert_p^p
$$
exists. Then clearly
$$
p^* \geq q \coloneqq \sup_{r>0} q(r) = \lim_{r\to\infty} q(r).
$$
We now distinguish two scenarios.

(i) If $q=0$, then $q(r)=0$ for every $r>0$. That is 
$$
\lim_{n\to\infty} \ \sup_{x\in\Gamma} \Vert {f_n}_{|B(x,r)} \Vert_p^p = 0 \quad \text{for every } r>0.
$$
Moreover, in this case we also have
$$
\lim_{n\to\infty} \Vert f_n^p \Vert_\infty = 0,
$$
which implies $\tf_n^p \to \tilde{0} \in \xtilde$:  indeed, using the representation \eqref{def-D}, the distance $\bD(\tf_n^p,\widetilde 0)$ between the orbits of the sub-probability measures $\tf_n^p$ and $\widetilde 0$ is bounded above by $\Vert f_n^p\Vert_1$ (recall that $\Vert f_n^p\Vert_1=\Vert f_n \Vert_p^p$ and thus $\Vert f_n^p\Vert_1\leq1$, so $\|f_n^p\|_1^k \leq \|f_n^p\|_1$ for $k\geq 1$), which goes to zero by the above display and the fact that $\sup_{n\in\N} \Vert f_n^p \Vert_1 \leq 1$. Hence, $\tf_n^p \to \tilde{0} \in \xtilde$. Thus $\tf_n \to \tilde{0} \in \xtilde_p$.

(ii) If $q>0$, then $q(r)>0$ for some $r>0$. We may thus find $r>0$ and $\{s_n\}_{n=1}^\infty \subset \Gamma$ such that
\begin{equation} \label{equ: lpdecomposition1}
\Vert \bigl(\delta_{s_n} * f_n \bigr)_{|B(e,r)} \Vert_p^p \geq \frac{q}{2} \quad \text{for all sufficiently large } n.
\end{equation}
Define $g_n \coloneqq \delta_{s_n}*f_n$. By going over to a subsequence we may assume that $g_n$ converges pointwise, as $n\to\infty$, to some $\alpha \in \Mcal_p^{\scriptscriptstyle{\leq 1}}$. It follows from Equation (\ref{equ: lpdecomposition1}) that
$$
\Vert \alpha\Vert_p^p \geq \Vert \alpha_{|B(e,r)} \Vert_p^p \geq \liminf_{n\to\infty} \Vert {g_n}_{|B(e,r)} \Vert_p^p \geq q/2.
$$
By a standard argument we may express $g_n = \alpha_n + \beta_n$ with $\alpha_n,\beta_n \in \Mcal_p^{\scriptscriptstyle{\leq 1}}$ such that
$$
\lim_{n\to\infty} \Vert \alpha_n - \alpha \Vert_p = 0
$$
and $\beta_n \to 0$ pointwise as $n\to\infty$. These properties imply that
\begin{flalign*}
p^*& =\lim_{n\to\infty} \Vert f_n \Vert_p^p = \lim_{n\to\infty} \Vert g_n \Vert_p^p = \lim_{n\to\infty} \sum_{x\in\Gamma} \bigl( \alpha_n(x) + \beta_n(x) \bigr)^p \\
& = \lim_{n\to\infty} \sum_{x\in B(e,r)} \bigl( \alpha_n(x) + \beta_n(x) \bigr)^p + \sum_{x\in B(e,r)^c} \bigl( \alpha_n(x) + \beta_n(x) \bigr)^p \\
& \geq \frac{q}{2} + \limsup_{n\to\infty} \Vert {\beta_n}_{|B(e,r)^c} \Vert_p^p \\
& = \frac{q}{2} + \limsup_{n\to\infty} \Vert \beta_n \Vert_p^p - \Vert {\beta_n}_{|B(e,r)} \Vert_p^p \\
& = \frac{q}{2} + \limsup_{n\to\infty} \Vert \beta_n \Vert_p^p
\end{flalign*}
and thus $\limsup_{n\to\infty} \Vert \beta_n \Vert_p^p \leq p^*-q/2$. Since $\alpha_n \to \alpha$ in $\ell^p(\Gamma)$, we may repeat the same computation for increasingly large balls around the identity to get the sharp estimate $\limsup_{n\to\infty} \Vert \beta_n \Vert_p^p \leq p^*-\Vert \alpha \Vert_p^p$. In any case, a positive amount of mass is removed when going over to the sequence $(\beta_n)_n$, which is the important point for the coming iterative procedure. In particular
$$
\limsup_{n\to\infty} q_{\beta_n}(r) \leq p^*-q/2 \quad \text{for every } r>0.
$$
Since $0\leq\beta_n \leq f_n$, we also have 
$$
\lim_{n\to\infty} q_{\beta_n}(r) \leq \lim_{n\to\infty} q_n(r)=q(r)\leq q \quad \text{for every } r>0
$$
and, putting the two estimates together, therefore
$$
\lim_{n\to\infty} q_{\beta_n}(r) \leq \min \{ q,p^*-q/2\} \quad \text{for every } r>0.
$$
We now shift this decomposition back by $\delta_{s_n^{-1}}$, i.e.\ take $\alpha_n^{(1)}=\delta_{s_n^{-1}} * \alpha_n$ and $\beta_n^{(1)}=\delta_{s_n^{-1}} * \beta_n$ giving a decomposition $f_n=\alpha_n^{(1)}+\beta_n^{(1)}$ such that
$$
\lim_{n\to\infty} \Vert \alpha - \delta_{s_n} * \alpha_n^{(1)} \Vert_p = 0
$$
and $\beta_n^{(1)}$ satisfies
$$
\lim_{n\to\infty} q_{\beta_n^{(1)}}(r) \leq \min \{ q,p^*-q/2\} \quad \text{for every } r>0
$$
and, since $\beta_n \hookrightarrow 0$, we also have
\begin{equation} \label{equ: prooflpdecomp1}
\lim_{n\to\infty} \Vert {\beta_n^{(1)}}_{|B(s_n,r)} \Vert_p^p = 0 \quad \text{for every } r>0.
\end{equation}
We may then apply the same procedure to the sequence $(\beta_n^{(1)})_n$. Note that if it is again possible to extract a positive amount of mass by applying suitable shifts, then Equation (\ref{equ: prooflpdecomp1}) ensures that these new locations of mass have to be widely separated from the locations $(s_n)_n$. Now iterate this procedure. If we are careful to remove at every step those pieces with the largest amount of mass, and finally pass to the diagonal sequence, we arrive at a subsequence with the desired properties, which we state as 
\begin{lemma}\label{lemma-decomposition}  Let $\Gamma$ be a countable discrete group and $p\in[1,\infty)$. Let $(f_n)_{n=1}^\infty$ be any sequence in $\Mcal_p$. Then there exists a subsequence such that the renamed sequence satisfies the following: 
There exist $\{\alpha_i,\alpha_i^{(n)}\}_{i,n\in\N} \subset \Mcal_p^{\scriptscriptstyle{\leq 1}}$ and $\{s_n^{(i)}\}_{i,n\in\N} \subset \Gamma$ and given any $\eps>0$ there exists $k\geq1$ such that one has the decomposition
$$
f_n = \sum_{i=1}^k \alpha_n^{(i)} + \gamma_{n}^{(k)}
$$
where $\{\gamma_{n}^{(k)}\}_{n\in\N} \subset \Mcal_p^{\scriptscriptstyle{\leq 1}}$, with the following properties:
\begin{enumerate}
\item The sequence $(\Vert \alpha_i \Vert_p)_i$ is non-increasing and $\sum_{i=1}^\infty \Vert \alpha_i \Vert_p^p \leq 1$.
\item For each $i \neq j$ we have that $(s_n^{(i)})^{-1}s_n^{(j)} \to \infty$ as $n \to \infty$.
\item For every $i\in\N$
$$
\lim_{n\to\infty} \Vert \delta_{s_n^{(i)}} * \alpha_n^{(i)} - \alpha_i \Vert_p = 0.
$$ 
\item The error $\gamma_n^{(k)}$ satisfies
$$
\limsup_{n\to\infty} \ \sup_{x\in\Gamma} \Vert {\gamma_n^{(k)}}_{|B(x,r)} \Vert_p^p \leq \eps \quad \text{for every } r>0
$$
and
$$
\lim_{n\to\infty} \Vert {\gamma_n^{(k)}}_{|B(s_n^{(i)},r)} \Vert_p^p = 0 \quad \text{for every } r>0 \text{ and } i \leq k,
$$
as well as
$$
\limsup_{n\to\infty} \Vert \gamma_n^{(k)} \Vert_p^p \leq 1-\sum_{i=1}^k \Vert \alpha_i \Vert_p^p,
$$
\end{enumerate}
\end{lemma}
To conclude the proof of Theorem \ref{theorem-CompactificationMp}, it remains to show that the decomposed subsequence converges in $\xtilde_p$. For this purpose, we will show that the sequence $\mu_n \coloneqq f_n^p$ satisfies $\tmu_n \to \xi = \eta^p$ in $\xtilde$. To see this we, define
$$
w_n^{(k)} = \mu_n - \sum_{i=1}^k \bigl( \alpha_n^{(i)} \bigr)^p = \Bigl( \sum_{i=1}^k {\alpha_n^{(i)}} + \gamma_n^{(k)} \Bigr)^p - \sum_{i=1}^k \bigl( \alpha_n^{(i)} \bigr)^p.
$$
Note that $w_n^{(k)} \in \Mcal_{\leq1}$ because $p\geq1$ and thus $\sum_{i=1}^k \bigl( \alpha_n^{(i)} \bigr)^p \leq \Bigl( \sum_{i=1}^k \alpha_n^{(i)} \Bigr)^p \leq f_n^p = \mu_n$. We then have the following properties:
\begin{enumerate}
\item The sequence $(\alpha_i^p(\Gamma))_i = (\Vert \alpha_i \Vert_p^p)_i$ is non-increasing and $\sum_{i=1}^\infty \alpha_i^p(\Gamma) \leq 1$.
\item For each $i \neq j$ we have that $(s_n^{(i)})^{-1}s_n^{(j)} \to \infty$ as $n\to\infty$.
\item For each $i\in\N$ we have
$$
\delta_{s_n^{(i)}} * \bigl(\alpha_n^{(i)}\bigr)^p \Rightarrow \alpha_i^p \quad \text{as } n \to \infty
$$
because $\delta_{s_n^{(i)}} * \alpha_n^{(i)} \to \alpha_i$ in $\ell^p(\Gamma)$.
\item By property (iii) above, each component $\alpha_n^{(i)}$ is concentrated on a ball around $s_n^{(i)}$. Since these sequences are widely separated, compare property (ii) above, and since we have only finitely many components we obtain that
$$
\lim_{n\to\infty} \Vert w_n^{(k)}- \bigl(\gamma_n^{(k)}\bigr)^p \Vert_\infty = 0.
$$
This immediately yields that for every fixed radius $r>0$
\begin{flalign*}
\limsup_{n\to\infty} \ \sup_{x\in\Gamma} \, w_n^{(k)}(B(x,r)) & = \limsup_{n\to\infty} \ \sup_{x\in\Gamma} \, \Vert {\bigl(\gamma_n^{(k)}\bigr)^p}_{|B(x,r)} \Vert_1 \\
& = \limsup_{n\to\infty} \ \sup_{x\in\Gamma} \, \Vert {\gamma_n^{(k)}}_{|B(x,r)} \Vert_p^p \leq \eps
\end{flalign*}
and for $i\leq k$ also
$$
\lim_{n\to\infty} \, w_n^{(k)}(B(s_n^{(i)},r)) = \lim_{n\to\infty} \Vert {\gamma_n^{(k)}}_{|B(s_n^{(i)},r)} \Vert_p^p = 0.
$$
\end{enumerate} 
From properties (i)-(iv) above it follows as in \cite[Proof of Theorem 3.2]{MV16} that the sequence $\mu_n$ convergences to the collection $\{\talpha_i^p\}=\eta^p$ in $\xtilde$, which concludes the proof. \eproof

It will be useful to consider a simplified version of Lemma \ref{lemma-decomposition} with fixed, but moving, components instead of converging components. To avoid summability issues, it is more convenient to simply consider approximations by finite sums.
\begin{cor} \label{cor: lpdecomposition} Let $p\in[1,\infty)$ and let $(f_n)_{n=1}^\infty$ be any sequence in $\Mcal_p$. Then there exists a subsequence such that the renamed sequence satisfies the following: 
There exist $\{\alpha_i\}_{i\in\N} \subset \Mcal_p^{\scriptscriptstyle{\leq 1}}$ and $\{s_n^{(i)}\}_{i,n\in\N} \subset \Gamma$ such that for every $\eps>0$ there exists $k\geq1$ such that one has the decomposition
$$
f_n = \sum_{i=1}^k \bigl( \delta_{s_n^{(i)}} * \alpha_i \bigr) + w_n^{(k)}
$$
where $\{w_{n}^{(k)}\}_{n\in\N} \subset \ell^p(\Gamma)$, with the following properties:
\begin{enumerate}
\item The sequence $(\Vert \alpha_i \Vert_p)_i$ is non-increasing and $\sum_{i=1}^\infty \Vert \alpha_i \Vert_p^p \leq 1$.
\item For each $i \neq j$ we have that $(s_n^{(i)})^{-1}s_n^{(j)} \to \infty$ as $n \to \infty$.
\item The error $w_n^{(k)}$ satisfies
$$
\limsup_{n\to\infty} \Vert w_n^{(k)} \Vert_\infty \leq \eps
$$
and
$$
\lim_{n\to\infty} \Vert {w_n^{(k)}}_{|B(s_n^{(i)},r)} \Vert_p^p = 0 \quad \text{for every } r>0 \text{ and } i \leq k,
$$
as well as
$$
\limsup_{n\to\infty} \Vert w_n^{(k)} \Vert_p^p \leq 1-\sum_{i=1}^k \Vert \alpha_i \Vert_p^p.
$$
\end{enumerate}
Moreover, in this case $\tf_n \to \{\talpha_i \}_{i\in I} \in \xtilde_p$, where $I=\{i \colon \alpha_i \neq 0\}$. 
\end{cor}

\begin{remark}
We point out that in the setting of Corollary \ref{cor: lpdecomposition} in particular
$$
\sup_{k\geq1} \ \limsup_{n\to\infty} \, \Vert w_n^{(k)} \Vert_p^p \leq 1.
$$
Applying log-convexity of $\ell^p$-norms, it follows that for every $\eps>0$ and every $q\in(p,\infty)$, there exists $k\geq1$ such that the error satisfies 
$$
\limsup_{n\to\infty} \Vert w_n^{(k)} \Vert_q \leq \eps.
$$
\end{remark}
\begin{remark} \label{remark: decompositionM1} We also record a simplified version of Lemma \ref{lemma-decomposition} for probability densities, which will be more convenient to use later.$^{\rm c}$\footnote{$^{\rm c}$  There are two main differences, which we want to point out: On the one hand, we may pass to an infinite sum as there are no summability issues when considering measures. On the other hand, we may again pass from convergence in each component to having fixed components which are being moved around by the group action. These two modifications may be allowed by absorbing the difference into an updated error term, which will in general not be a measure anymore as it can have negative values. However, the error terms will still form a uniformly bounded sequence of integrable functions which vanish uniformly over the group. In this case the limiting $\ell^1$-norm can be easily derived.}
Let $(\mu_n)_{n=1}^\infty$ be any sequence in $\Mcal_1(\Gamma)$. Then there exists a subsequence such that for the renamed sequence one has the decomposition
\begin{equation} \label{equ: discretedecomposition}
\mu_n = \sum_{i=1}^\infty \delta_{s_n^{(i)}} * \alpha_{i} + w_n,
\end{equation}
where $\{ \alpha_i\}_{i \in \N} \subset \Mcal_1^{\scriptscriptstyle{\leq 1}}(\Gamma)$, $\{s_n^{(i)}\}_{i,n\in\N} \subset \Gamma$ and $\{w_n\}_{n\in\N} \subset \ell^1(\Gamma)$, satisfying the following properties: 
\begin{enumerate}
\item The sequence $(\alpha_i(\Gamma))_i$ is non-increasing and $\sum_{i=1}^\infty \alpha_i(\Gamma) \leq 1$.
\item For each $i \neq j$, we have that $(s_n^{(i)})^{-1}(s_n^{(j)}) \to \infty$ as $n \to \infty$, in the sense that the sequence escapes every finite set.
\item The error $w_n$ converges to zero in $\ell^\infty(\Gamma)$ as $n \to \infty$ and satisfies
$$
\limsup_{n\to\infty} \Vert w_n \Vert_1 = 1-\sum_{i=1}^\infty \alpha_i(\Gamma).   				
$$
\end{enumerate}
\end{remark}

\begin{remark} Instead of verifying the existence of the decompositions in Lemma \ref{lemma-decomposition} in the way presented above, we can also deduce it by first showing the simpler version of this statement 
for $p=1$ and then, by considering the sequence $(f_n^p)_n$ instead, ``pull back‘‘ the decomposition by passing to the $p^{th}$-root of each component. We chose the above approach because it seems conceptually clearer. 
\end{remark}
\begin{remark} 
The case $p=2$ is particularly interesting because of the ambient Hilbert space $\ell^2(\Gamma)$. In this setting, it is not difficult to verify the  properties of the decomposition in Corollary \ref{cor: lpdecomposition}
because of the availability of Pythagorian identity. 
Building on an idea using Levy's concentration function, 
decompositions of measures of this type date back to a seminal work of Parthasarath{\color{blue}y}, Rao and Varadhan \cite{PRV62} 
which have been developed further. 
It seems to us that it remained somewhat elusive whether such decompositions may be understood as a special case of classical compactness. We believe that the group-invariant compactification introduced in \cite{MV16} provides a satisfactory answer to this question (at least for the purposes we have in mind) -- namely, it identifies the quotient 
as the right space to be compactified and constructs $\widetilde{\mathcal X}$ with a topology so that 
the quotient can be compactly embedded into $\xtilde$. 
There is a price to pay, of course, which is that (by modding out by the group action) we cannot retrieve the exact locations of the components, neither the 
precise form of the mass that totally dissipates. However, these ``lost" information do not seem to impact the relevant purposes and the applications, as explored also in this article in the context of $C^*$-algebras and geometric properties of groups. 
\end{remark}

\section{Asymptotically orthogonal decompositions of Schur multipliers}\label{main}

\subsection{The first main result.}

This section is devoted to our first main result, Theorem \ref{theorem-MainTheorem}, which relates the compactification of $\tMcal_p$ to decompositions of certain Schur multipliers on the uniform Roe algebra and the reduced group $C^*$-algebra of a countable discrete group $\Gamma$. Let us first start with the precise definition of {\it asymptotically orthogonal decomposition} of a sequence of bounded linear operators $(T_n)_{n=1}^\infty \subset B(X,Y)$ w.r.t. a given collection $\mathcal S\subset B(X,Y)$:

\begin{definition}\label{def-orthogonal}
For Banach spaces $X$ and $Y$, let $B(X,Y)$ denote the Banach space of bounded linear operators from $X$ to $Y$. The {\it strong operator topology} on $B(X,Y)$ is the topology induced by the seminorms 
$$
T \mapsto \Vert Tx \Vert, \quad \mbox{for} \, \, x \in X.
$$
Given a sequence $(T_n)_{n=1}^\infty$ in $B(X,Y)$ and a subset $\skris \subset B(X,Y)$, an {\it asymptotically orthogonal decomposition} of $(T_n)_n$ with respect to $\skris$ is a norm-convergent sum
$$
\sum_{i \in I} S_i,
$$
with $I \subset \N$ (the empty sum is identified with the zero operator), which has the following properties:
\begin{enumerate}
\item $S_i \in \skris$ for every $i \in I$.
\item There exists a subsequence of $(T_n)_n$ such that the renamed sequence satisfies
$$
\lim_{n\to\infty} T_n = \sum_{i\in I} S_i
$$
in the strong operator topology on $B(X,Y)$ and such that the following asymptotic norm equality holds
$$
\lim_{n\to\infty} \Vert T_n \Vert = \sum_{i\in I} \Vert S_i \Vert.
$$
\end{enumerate} 
\end{definition}

Recall the definition of the positive definite functions $\Phi(f)$ and that of the corresponding Schur multipliers $M_{\Phi(f)}$ from Section \ref{sec-Schur}. In the sequel we will also write
\begin{equation}\label{Mxi}
M_\xi= \sum_{\widetilde\alpha\in \xi} M_{\Phi(\alpha)}, \qquad \mbox{for} \, \, \xi=\{\widetilde\alpha_i\}_{i\in I}\in \widetilde{\mathcal X}_p.
\end{equation}

Here is our first theorem. 
\begin{theorem} \label{theorem-MainTheorem}
Let $\Gamma$ be a countable discrete group and $p \in [1,2)$. For any sequence $(f_n)_{n=1}^\infty$ in $\Mcal_p(\Gamma)$, consider the sequence $(M_n)_{n=1}^\infty$ of Schur multipliers $M_n \coloneqq M_{\Phi(f_n)}$ on $C_u^*(\Gamma)$
or on $C_\lambda^*(\Gamma)$. Then the set of asymptotically orthogonal decompositions of $(M_n)_n$ with respect to 
$$
\left\{ M_{\Phi(f)} : f \in \Mcal_p^{\scriptscriptstyle{\leq 1}}(\Gamma) \right\}
$$
is given by 
$$
\left\{ M_\xi \colon \xi \text{ is a limit point of } (\tf_n)_n \text{ in } \xtilde_p(\Gamma) \right\}
\vspace{2mm}
$$ 
and contains all limit points of the sequence $(M_n)_n$ in the strong operator topology. Moreover, whenever
$$
\lim_{n \to \infty} M_{n} = M
$$
exists in the strong operator topology, every limit point $\xi \in \xtilde_p(\Gamma)$ of $(\tf_{n})_n$ satisfies $M_\xi= M$.
\end{theorem}

Combining Theorem \ref{theorem-MainTheorem} and compactness of $\xtilde_p(\Gamma)$ we immediately obtain the following corollary.

\begin{cor}\label{corollary-compactnessM1} Let $\Gamma$ be a countable discrete group and $p\in [1,2)$. Then 
$$
\bS_p(\Gamma) = \left\{ M_{\Phi(f)} \colon f \in \Mcal_p(\Gamma) \right\} \subset B(C^*_u(\Gamma))
$$ 
is relatively sequentially compact with respect to the strong operator topology and its closure in the strong operator topology is given by 
$$
\bX_p(\Gamma)=\{ M_\xi \colon \xi \in \xtilde_p(\Gamma) \}.
$$
\end{cor}

Regarding norm convergence, some adaptations to our method of proving Theorem \ref{theorem-MainTheorem} yield the following characterization.

\begin{theorem}\label{thm-normconvergence} Let $\Gamma$ be a countable discrete group and $p\in[1,2)$. For any sequence $(f_n)_{n=1}^\infty$ in $\Mcal_p(\Gamma)$, consider the sequence $(M_n)_{n=1}^\infty$ of Schur multipliers $M_n\coloneqq M_{\Phi(f_n)}$ on $B(\ell^2(\Gamma))$. Then $(M_n)_n$ converges in norm if and only if there exists $\alpha\in\Mcal_p^{\scriptscriptstyle{\leq 1}}(\Gamma)$ such that
$$
\lim_{n\to\infty} \widetilde f_n = \{\widetilde \alpha\} \quad \mbox{in} \ \ \xtilde_p(\Gamma).
$$
\end{theorem}

As mentioned in the introduction, let us also record the following consequences of Theorem \ref{theorem-MainTheorem} regarding positive definite functions.

\begin{cor}\label{corollary-limitpoints} Let $\Gamma$ be a countable discrete group and $p\in[1,2)$. Then any pointwise limit $\varphi$ of a sequence in $P_p(\Gamma)$ is of the form
$$
\varphi = \sum_{i=1}^\infty \Phi(f_i), \qquad\mbox{where}\quad f_i \in \Mcal_p^{\scriptscriptstyle{\leq 1}},\,\,\mbox{with}\,\, \sum_{i=1}^\infty \|f_i\|_p^p \leq 1.
$$
\end{cor}

\begin{cor}\label{corollary-minimizers}  Let $\Gamma$ be a countable discrete group and $p\in[1,2)$. Then for any finite subset $F\subset\Gamma$, there exist $f_i \in \Mcal_p^{\scriptscriptstyle{\leq 1}}, i \in \N$, with $\sum_{i=1}^\infty \|f_i\|_p^p \leq1$ such that the positive definite function
$$
\varphi=\sum_{i=1}^\infty \Phi(f_i)
$$
minimizes the expression
$$
d_F(\psi,1) \coloneqq \sup \bigl\{ |\psi(s)-1| \colon s \in F \bigr\}
$$
over all positive definite functions $\psi$ of this form and in particular over all $\psi \in P_p(\Gamma)$.
\end{cor}

In Section \ref{sec-proof-main} we will prove Theorem \ref{theorem-MainTheorem}, while the proofs of Theorem \ref{thm-normconvergence} and Corollary \ref{corollary-limitpoints}-\ref{corollary-minimizers} will be contained in Section \ref{sec-proof-appli}.

\subsection{Proof of Theorem \ref{theorem-MainTheorem}.}\label{sec-proof-main} 

Throughout this subsection let $\Gamma$ be a countable discrete group. To shorten notation, the dependence on the group will again be dropped in our notations, so e.g.\ $\Mcal_1=\Mcal_1(\Gamma)$. The proof proceeds in the following way: We first prove a convergence criterion for Schur multipliers and then verify the criterion in the set-up of our main theorem. Finally we compute the norms of our Schur multipliers precisely to prove the asymptotic norm equality.

\subsubsection{Strong convergence of Schur multipliers.} Recall the following basic result.

\begin{lemma}[Uniform Boundedness Principle] \label{lemma: UniformBoundednessPrinciple}
Let $X$ be a Banach space, let $D \subset X$ a dense subspace and let $(T_n)_{n=1}^\infty$ and $T$ be operators in $B(X)$. Then the following are equivalent:
\begin{enumerate}
\item The sequence $(T_n)_n$ converges to $T$ in the strong operator topology on $B(X)$.
\item The sequence $(\Vert T_n \Vert)_n$ is bounded and the restriction of $T_n$ to $D$ converges to the restriction of $T$ to $D$ in the strong operator topology on $B(D,X)$.
\end{enumerate}
\end{lemma}

From the Uniform Boundedness Principle we can derive the following characterization of convergence of Schur multipliers in the strong operator topology.

\begin{lemma} \label{lemma-strongconvergence}
Let $(\varphi_n)_{n=1}^\infty$ and $\varphi$ be functions on $\Gamma$ such that $(M_{\varphi_n})_{n=1}^\infty$ and $M_\varphi$ are bounded operators on $C^*_u(\Gamma)$. Then the following are equivalent:
\begin{enumerate}
\item The sequence $(M_{\varphi_n})_n$ converges to $M_{\varphi}$ in the strong operator topology on $B(C^*_u(\Gamma))$.
\item The sequence $(\Vert M_{\varphi_n} \Vert)_n$ is bounded and $\varphi_n$ converges pointwise to $\varphi$.
\end{enumerate}
\end{lemma}

\proof To shorten notation, define $M_n \coloneqq M_{\varphi_n}$. Recall that $C_u^*(\Gamma)$ is the norm closure of the $^*$-algebra $\mathbb C_u[\Gamma]$ of operators with finite propagation. By Lemma \ref{lemma: UniformBoundednessPrinciple} it therefore suffices to prove that if
$$
\sup_{n \in \N} \Vert M_n \Vert < \infty,
$$
then the following are equivalent:
\begin{enumerate}[label=(\Roman*)]
\item The sequence of restrictions of $M_n$ to $\mathbb C_u[\Gamma]$ converges to the restriction of $M_{\varphi}$ to $\mathbb C_u[\Gamma]$ in the strong operator topology on $B(\mathbb C_u[\Gamma])$.
\item The sequence $(\varphi_n)_n$ converges pointwise to $\varphi$.
\end{enumerate}

To prove the first implication, fix some $x \in \Gamma$ and note that for every function $\psi$ on $\Gamma$ we have $\psi(x)=K_\psi(x,e)$. Now define
$$
T(s,t) = \1_{(x,e)}(s,t).
$$
Clearly $T \in \mathbb C_u[\Gamma]$ and thus $M_n(T) \to M_\varphi(T)$ in $B(\ell^2(\Gamma))$. In matrix form
$$
M_n(T) = [\1_{(x,e)}(s,t) \varphi_n(st^{-1})]_{(s,t) \in \Gamma \times \Gamma} \quad \text{and} \quad M_\varphi(T) = [\1_{(x,e)}(s,t) \varphi(st^{-1})]_{(s,t) \in \Gamma \times \Gamma},
$$
hence it follows from the Cauchy Schwarz inequality that
$$
|\varphi_n(x)-\varphi(x)| = \bigl| \langle \delta_x, M_n(T)\delta_e \rangle - \langle \delta_x, M_\varphi(T)\delta_e \rangle \bigr| \leq \Vert M_n(T) - M_\varphi(T) \Vert,
$$
which goes to zero as $n \to \infty$. This proves pointwise convergence of $(\varphi_n)_n$ to $\varphi$. 

To prove the reverse implication, let $T \in \mathbb C_u[\Gamma]\setminus\{0\}$ and $\eps>0$ be arbitrary. We will verify that 
\begin{equation} \label{claim-proofconvlemma1}
\lim_{n \to \infty} \Vert M_n(T) - M_\varphi(T) \Vert = 0.
\end{equation} 
Choose a finite subset $F \subset \Gamma$ such that $\supp(T) \subset \Tube(F)$, then
$$
T(x,y)=0 \quad \text{ for all } (x,y) \in \Gamma \times \Gamma \text{ such that } xy^{-1} \in \Gamma \setminus F.
$$
Hence also for every function $\psi \in \ell^\infty(\Gamma \times \Gamma)$
$$
\supp(\psi \circ T) \subset \Tube(F).
$$
Note that if $\psi$ is of the form $K_f$ for some function $f$, then the values of $K_f$ on $\Tube(F)$ depend only on the values of $f$ on the finite set $F$. Since $\varphi_n$ converges pointwise to $\varphi$, the corresponding restrictions to $F$ converge uniformly. Therefore
$$
\sup_{x \in F} | \varphi_n(x) - \varphi(x) | < \frac{\eps}{\Vert T \Vert_\infty|F| } \quad \text{for all sufficiently large } n.
$$
This means precisely that for all sufficiently large $n$, the operator $M_n(T) - M_\varphi(T)$, which has the matrix representation
$$
\bigl[\bigl(\varphi_n(xy^{-1})-\varphi(xy^{-1})\bigr)T(x,y)\bigr]_{(x,y) \in \Gamma \times \Gamma},
$$
is supported in $\Tube(F)$ and satisfies $\Vert M_n(T) - M_\varphi(T) \Vert_\infty < \eps/|F|$. Applying Lemma \ref{lemma: FinPropBound} shows that
$$
\Vert M_n(T) - M_\varphi(T) \Vert \leq \eps
$$
for all sufficiently large $n$ which proves (\ref{claim-proofconvlemma1}). Since $T \in \mathbb C_u[\Gamma]\setminus\{0\}$ was arbitrary $M_n$ converges to $M_\varphi$ in the strong operator topology on $B(\mathbb C_u[\Gamma])$ and the proof of Lemma \ref{lemma-strongconvergence} is complete.
\eproof

\subsubsection{Continuity of the map $\Phi$.} Recall that the construction of passing from $f\in\Mcal_1$ to $M_{\Phi(f)}$ depended only on the orbit of $f$ under the left action of $\Gamma$. The following lemma is an important
 ingredient in the proof of Theorem \ref{theorem-MainTheorem}. Roughly speaking, it asserts that our construction is continuous on $\tMcal_1$.

\begin{lemma}\label{lemma-ContinuityPhi} Let $(f_n)_{n=1}^\infty$ be a sequence in $\Mcal_1$ such that 
$$
\lim_{n\to\infty} \tf_n = \xi = \{\talpha_i\}_{i\in I} \in \xtilde.
$$
Then 
$$
\lim_{n\to\infty} \Phi(f_n) = \sum_{i \in I} \Phi(\alpha_i) \quad \text{pointwise}.
$$
\end{lemma}

\proof Recall that by Lemma \ref{lemma-decomposition} and Remark \ref{remark: decompositionM1}, for any sequence $(\mu_n)_{n=1}^\infty$ in $\Mcal_1$ there exists a subsequence such that for the renamed sequence we have the decomposition
$$
\mu_n = \sum_{i=1}^\infty \delta_{t_n^{(i)}} * \beta_{i} + h_n,
$$
where $\{t_n^{(i)}\}_{i,n\in\N} \subset \Gamma$ are shifts, $\{ \beta_i\}_{i \in I} \subset \Mcal_{\leq 1}(\Gamma)$ are measures and $\{h_n\}_{n\in\N} \subset \ell^1(\Gamma)$ are functions satisfying the following properties: 
\begin{enumerate}
\item The sequence $(\beta_i(\Gamma))_i$ is non-increasing and $\sum_{i=1}^\infty \beta_i(\Gamma) \leq 1$.
\item For each $i \neq j$, we have that $(t_n^{(i)})^{-1}t_n^{(j)} \to \infty$ as $n \to \infty$.
\item The error $h_n$ converges to zero in $\ell^\infty(\Gamma)$ as $n \to \infty$ and satisfies
$$
\limsup_{n\to\infty} \Vert h_n \Vert_1 = 1-\sum_{i=1}^\infty \beta_i(\Gamma).
$$
\end{enumerate}
For a subsequence of $(f_n)_n$, the fact that $\tf_n$ converges to $\xi$ implies that we can match $\beta_i = \alpha_i$ in the above decomposition. Applying this argument to every subsequence, we deduce that there exist shifts $\{s_n^{(i)}\}_{i\in I,n\in\N} \subset \Gamma$ such that we have the decomposition
\begin{equation} \label{equ: LemmaContPhi1}
f_n = \sum_{i=1}^\infty \delta_{s_n^{(i)}} * \alpha_{i} + w_n,
\end{equation}
where $\{s_n^{(i)}\}_{i,n\in\N} \subset \Gamma$ are shifts and $\{w_n\}_{n\in\N} \subset \ell^1(\Gamma)$ are functions satisfying the following properties: 
\begin{enumerate}
\item For each $i \neq j$, we have that $(s_n^{(i)})^{-1}s_n^{(j)}\to \infty$ as $n \to \infty$.
\item The error $w_n$ converges to zero in $\ell^\infty(\Gamma)$ as $n \to \infty$ and satisfies
$$
\limsup_{n\to\infty} \Vert w_n \Vert_1 = 1-\sum_{i=1}^\infty \alpha_i(\Gamma).
$$
\end{enumerate}
Recall that we can write 
$$
\Phi(f_n)(x) = \langle f_n, \rho(x)f_n \rangle \quad \text{for every } x\in\Gamma.
$$
Now fix $x\in\Gamma$. Setting 
$$
\alpha_n^{(i)} \coloneqq \delta_{s_n^{(i)}} * \alpha_i
$$ 
and plugging the decomposition (\ref{equ: LemmaContPhi1}) into the inner product form, we obtain that
\begin{equation} \label{equ: LemmaContPhi2}
\Phi(f_n)(x) = \langle w_n, \rho(x)f_n \rangle + \langle f_n - w_n, \rho(x)w_n \rangle + \sum_{i \in I} \langle \alpha_n^{(i)}, \rho(x)\alpha_n^{(i)} \rangle + \sum_{i,j\in I,i\neq j} \langle \alpha_n^{(i)}, \rho(x)\alpha_n^{(j)} \rangle.
\end{equation}
Using invariance of $\Phi$ under left translations we see that
$$
\sum_{i \in I} \langle \alpha_n^{(i)}, \rho(x)\alpha_n^{(i)} \rangle = \sum_{i\in I} \Phi(\alpha_n^{(i)})(x) = \sum_{i\in I} \Phi(\alpha_i)(x).
$$
Therefore it remains to show that the other terms vanish in the limit $n\to\infty$. 

Since the error is uniformly bounded in $\ell^1(\Gamma)$ and $\Vert w_n \Vert_\infty \to 0$, it follows from log-convexity of $\ell^p$-norms that $\Vert w_n \Vert_p \to 0$ for every $1<p\leq \infty$. In particular it converges to zero in $\ell^2(\Gamma)$, which proves that the first two summands in Equation (\ref{equ: LemmaContPhi2}) vanish in the limit $n\to \infty$.

Regarding the last term in Equation (\ref{equ: LemmaContPhi2}), we first note that $\sum_{i\in I} \Vert \alpha_n^{(i)} \Vert_1 \leq 1$ implies that for every $\eps>0$ we find a finite subset $J$ of $I$ such that
$$
\sum_{i\in I\setminus J} \Vert \alpha_n^{(i)} \Vert_1 \leq \eps.
$$
Then in particular
\begin{flalign*}
\sum_{i,j \in I \setminus J} \langle \alpha_n^{(i)}, \rho(x)\alpha_n^{(j)} \rangle & = \sum_{i,j \in I \setminus J} \ \sum_{y\in\Gamma}  \alpha_n^{(i)}(y)\rho(x)\alpha_n^{(j)}(y) \leq \sum_{i,j \in I \setminus J} \Vert \rho(x)\alpha_n^{(j)} \Vert_\infty \sum_{y\in\Gamma} \alpha_n^{(i)}(y) \\
& \leq \sum_{i,j \in I \setminus J} \Vert \alpha_j \Vert_1 \Vert \alpha_i \Vert_1 \leq  \eps^2,
\end{flalign*}
where the final estimate follows by interchanging the order of summation which is possible because all summands are non-negative. Hence it suffices to prove that for every finite $J$
$$
\lim_{n\to\infty} \sum_{i,j\in J,i\neq j} \langle \alpha_n^{(i)}, \rho(x)\alpha_n^{(j)} \rangle = 0.
$$
Since the sum is finite, it is enough to check that every summand converges to zero. Note that the measure $\alpha_n^{(i)}$ is concentrated on a large ball of fixed finite radius $r_i>0$ near $s_n^{(i)}$, while $\rho(x)\alpha_n^{(j)}$ is concentrated on a ball of fixed finite radius $r_j>0$ near $s_n^{(j)}x^{-1}$. By property (ii) above, we also have that $(s_n^{(i)})^{-1}s_n^{(j)}x^{-1} \to \infty$ and thus the two balls $B_{i,n} \coloneqq B(s_n^{(i)},r_i)$ and $B_{j,n} \coloneqq B(s_n^{(j)}x^{-1},r_j)$ are disjoint for all large enough $n$. Now consider the rough estimate
$$
\langle \alpha_n^{(i)}, \rho(x)\alpha_n^{(j)} \rangle = \sum_{y\in\Gamma}  \alpha_n^{(i)}(y)\rho(x)\alpha_n^{(j)}(y) \leq \Vert (\rho(x)\alpha_n^{(j)})_{|B_{i,n}} \Vert_\infty + \alpha_n^{(i)}(B_{i,n}^c).
$$
By the argument above the upper bound can be made arbitrarily small for large enough $n$ by choosing a suitably large ball, whence $\lim_{n\to\infty} \langle \alpha_n^{(i)}, \rho(x)\alpha_n^{(j)} \rangle = 0$. This concludes the proof of Lemma \ref{lemma-ContinuityPhi}. \eproof

\subsubsection{Computing the norms.} For future use, we look at a slightly more general set-up than what is needed to prove Theorem \ref{theorem-MainTheorem}. To that end, recall that the map $\Phi$ was defined on $\ell^2$-functions via
\vspace{1mm}
$$
\Phi(f)(x) \coloneqq \sum_{y \in \Gamma} f(y)\overline{f(yx)} = \langle f,\rho(x)f \rangle \quad \text{for } f \in \ell^2(\Gamma).
$$

\begin{lemma} \label{lemma: PhifPosDef} Let $ f \in \ell^2(\Gamma)$, then $\Phi(f)$ is positive definite. Moreover, if $f$ is finitely supported, then so is $\Phi(f)$.
\end{lemma}

\proof Positive definiteness follows because we can rewrite $\Phi(f)$ as
$$
\Phi(f)(x) = \langle f, \rho(x) f \rangle,
$$
i.e.\ the matrix coefficient of a unitary representation, see e.g.\ \cite[Theorem 2.5.11]{BrownOzawa}. Now let $f\in C_c(\Gamma)$ and choose $R \geq 0$ such that $\supp f \subset B(e,R)$. We claim that $\supp \, \Phi(f) \subset B(e,2R)$, which implies that $\Phi(f)$ has finite support by properness of $d_\Gamma$. To see this let $x\in\Gamma$ such that $d_\Gamma(e,x)>2R$. If $y\in\Gamma$ is such that $f(y)\neq0$, then $d_\Gamma(e,y)\leq R$. In this case, the reverse triangle inequality and left-invariance imply that
$$
d_\Gamma(e,yx) \geq \bigl| d_\Gamma(e,y) - d_\Gamma(y,yx) \bigr| = \bigl| d_\Gamma(e,y) - d_\Gamma(e,x) \bigr| > R
$$
and thus $f(yx)=0$. Therefore $\Phi(x)=\sum_{y \in \Gamma} f(y)\overline{f(yx)}=0$, which proves the claim. \eproof.

\begin{cor} \label{cor: HaagerupLemma} Let $f \in \ell^2(\Gamma)$, then 
$$
M_{\Phi(f)} \colon B(\ell^2(\Gamma)) \longrightarrow B(\ell^2(\Gamma))
$$
is completely positive with norm $\Vert M_{\Phi(f)} \Vert = \Vert f \Vert_2^2$. \\
Moreover, its restrictions to $C_u^*(\Gamma)$ and $C_\lambda^*(\Gamma)$ are both completely positive with norm equal to $\Vert f \Vert_2^2$.
\end{cor}

\proof Since
$$
\Phi(f)(e) = \sum_{y \in \Gamma} f(y)\overline{f(y)} = \Vert f \Vert_2^2,
$$
the claim follows from Lemma \ref{lemma: PhifPosDef} and Lemma \ref{lemma: GeneralizedHaagerup} . \eproof

\subsubsection{Concluding the proof.} We now combine the previous steps to prove Theorem \ref{theorem-MainTheorem}. Let us first handle the case $p=1$ for simplicity. 

Let $(f_n)_{n=1}^\infty$ be any sequence in $\Mcal_1(\Gamma)$ and consider the sequence $(M_n)_{n=1}^\infty$ of Schur multipliers $M_n \coloneqq M_{\Phi(f_n)}$ on $C_u^*(\Gamma)$. Suppose $\xi=\{\talpha_i\}_{i\in I}$ is a limit point of $(\tf_n)_n$ in $\xtilde$, then Lemma \ref{lemma-ContinuityPhi} implies that, for some subsequence, we have $\lim_{k\to\infty} \Phi(f_{n_k}) = \sum_{i \in I} \Phi(\alpha_i)$ pointwise. Since the operators $M_\xi$ and $M_n$ are uniformly bounded, Lemma \ref{lemma-strongconvergence} implies that 
$$
\lim_{k\to\infty} M_{n_k} = M_\xi
$$
in the strong operator topology on $B(C_u^*(\Gamma))$. Finally, using the fact that sums of positive definite functions are positive definite functions and applying Corollary \ref{cor: HaagerupLemma} twice, we obtain
$$
\lim_{n\to\infty} \Vert M_{n_k} \Vert = \lim_{n\to\infty} \Phi(f_{n_k})(e) = \sum_{i\in I} \Phi(\alpha_i)(e) = \sum_{i\in I} \Vert M_{\Phi(\alpha_i)} \Vert = \Vert M_\xi \Vert.
$$
That is, we have shown that $M_\xi$ is an asymptotically orthogonal decomposition of $(M_n)_n$. 

On the other hand, suppose that $M$ is any strong operator topology limit point of $(M_n)_n$, say along the subsequence $(M_{n_k})_k$. By compactness of $\xtilde$, the sequence $(\tf_{n_k})_k$ has another subsequence which converges to some $\xi \in \xtilde$. Our previous argument implies that along this sub-subsequence the corresponding Schur multipliers converge to $M_\xi$ and hence $M=M_\xi$. Therefore every strong operator topology limit point of the sequence $(M_n)_n$ can be written as such an asymptotically orthogonal decomposition.

Let us now handle the case $p\in (1,2)$, for which we will follow a similar strategy as above and only indicate the necessary changes: Lemma \ref{lemma-strongconvergence} and Corollary \ref{cor: HaagerupLemma}, as well as compactness of $\xtilde_p$, can be used in the same way, provided that continuity guaranteed by Lemma \ref{lemma-ContinuityPhi} is still valid. Thus essentially the only difference lies in checking that for any sequence $(f_n)_{n=1}^\infty$ in $\Mcal_p$ such that $\tf_n \to \{\talpha_i\}_{i\in I} \in \xtilde_p$, again
$$
\lim_{n \to \infty} \Phi(f_n) = \sum_{i\in I} \Phi(\alpha_i) \quad \text{pointwise}.
$$
Note that the right-hand side exists because the Cauchy-Schwarz inequality implies
$$
\sum_{i\in I} \Phi(\alpha_i)(x) = \sum_{i\in I} \langle \alpha_i, \rho(x)\alpha_i \rangle \leq \sum_{i\in I} \Vert \alpha_i \Vert_2^2 \leq \sum_{i\in I} \Vert \alpha_i \Vert_p^p \leq 1.
$$
Using now the decomposition guaranteed by Corollary \ref{cor: lpdecomposition} we obtain that for any $\eps>0$ there exists a finite set $J \subset I$ such that uniformly
$$
\Bigl| \sum_{i\in I} \Phi(\alpha_i) - \sum_{i\in J} \Phi(\alpha_i) \Bigr| < \eps
$$
and 
$$
f_n = \sum_{i \in J} \bigl( \delta_{s_n^{(i)}} * \alpha_i \bigr) + w_n^{(J)},
$$
where $\{s_n^{(i)}\}_{i\in J,n\in N}$ and $\{w_{n}^{(J)}\}_{n\in\N} \subset \ell^p(\Gamma)$, with the following properties:
\begin{enumerate}
\item For each $i \neq j$ we have that $(s_n^{(i)})^{-1}s_n^{(j)} \to \infty$ as $n \to \infty$.
\item The error satisfies $\limsup_{n\to\infty} \Vert w_n^{(J)} \Vert_2 \leq \eps$. Note that this is possible because $p < 2$.
\end{enumerate}
Now $\Phi(f_n)$ has the same expression as in Equation (\ref{equ: LemmaContPhi2}) and it follows from the Cauchy-Schwarz inequality that the terms featuring the error are negligible. Finally, the (finite) sum over the interaction terms between profiles converges to zero as $n\to\infty$ because of asymptotic orthogonality in $\ell^2(\Gamma)$. The proof of Theorem \ref{theorem-MainTheorem} is thus complete. \eproof

\begin{remark} To have an analogue of Theorem \ref{theorem-MainTheorem} beyond $p>2$, a different construction would be necessary due to the use of the inner product in the definition of the map $\Phi$. We will not pursue this.
\end{remark}

\begin{remark} \label{remark: CriticalCase}
Let us also comment on the critical case $p=2$. While our construction still makes sense, the analogous result is not valid anymore, which is essentially due to the fact that totally disintegrating sequences do not converge to zero in $\ell^2(\Gamma)$. In fact, for amenable groups, there are examples of such sequences such that the corresponding Schur multipliers converge to the identity. Fortunately, the compactifications of $\bS_p$ with $p\in[1,2)$ can approximate this case, see Section \ref{section-applications} for details.
\end{remark}

\subsection{Proofs of Theorem \ref{thm-normconvergence} and Corollary \ref{corollary-limitpoints}-\ref{corollary-minimizers}.}\label{sec-proof-appli}

We start with the proof of Theorem \ref{thm-normconvergence}.

\noindent{\bf Proof of Theorem \ref{thm-normconvergence}.} Suppose there exists $\alpha \in \Mcal_p^{\scriptscriptstyle{\leq1}}$ such that 
$$
\lim_{n\to\infty} \widetilde f_n = \{ \widetilde \alpha \}.
$$
Then as carried out above in the proof of Theorem \ref{theorem-MainTheorem}, we may decompose 
$$
f_n = \delta_{s_n} * \alpha + w_n
$$
where $(s_n)_{n=1}^\infty \subset \Gamma$ and $(w_n)_{n=1}^\infty \subset \ell^p(\Gamma)$ such that $\Vert w_n \Vert_2 \to 0$ as $n\to\infty$. Applying this decomposition, we obtain that
$$
\Phi(f_n)(s)=\Phi(\alpha)(s) + \langle \delta_{s_n}*\alpha, \rho(s) w_n \rangle + \langle w_n, \rho(s)(\delta_{s_n} * \alpha) \rangle + \Phi(w_n)(s).
$$
Define the functions
$$
h_n(s) \coloneqq \langle \delta_{s_n}*\alpha, \rho(s) w_n \rangle \quad \mbox{and} \quad h_n'(x) \coloneqq \langle w_n, \rho(s)(\delta_{s_n} * \alpha) \rangle
$$
and let $k_n(s,t)\coloneqq h_n(st^{-1})$ and $k_n'(s,t)\coloneqq h_n'(st^{-1})$ be the corresponding kernels. We can rewrite
\begin{equation}\label{Mn}
M_n = M_{\Phi(\alpha)} + M_{h_n} + M_{h_n'} + M_{\Phi(w_n)}
\end{equation}
and it suffices to prove that the latter three operators converge to zero in norm. Note that $\Phi(w_n)$ is positive definite, therefore by Lemma \ref{lemma: GeneralizedHaagerup} 
$$
\Vert M_{\Phi(w_n)} \Vert = \Phi(w_n)(e)=\Vert w_n \Vert_2^2 \to 0
$$
 as $n\to\infty$. Now let us consider the second and the third summand on the right hand side in \eqref{Mn}. Note that the functions $h_n$ and 
$h_n'$ are in general not positive definite. However, defining
$$
\varphi_t := \rho(t)(\delta_{s_n} * \alpha) \quad\mbox{and}\quad \psi_s := \rho(s)w_n 
$$
we have that
$$
k_n(s,t)=h_n(st^{-1})=\langle \rho(t)(\delta_{s_n}*\alpha), \rho(s) w_n \rangle \quad \mbox{for all} \ s,t \in \Gamma.
$$
By Theorem \ref{theorem: Grothendieck} it follows that
$$
\Vert M_{h_n} \Vert \leq \sup_{t\in\Gamma} \Vert \varphi_t \Vert_2 \ \sup_{s\in\Gamma} \Vert \psi_s \Vert_2 = \Vert \alpha \Vert_2 \Vert w_n \Vert_2,
$$
which converges to zero as $n\to\infty$. A similar argument shows that $\|M_{h_n'}\|\to 0$ as $n\to\infty$, which proves that $M_n$ converges in norm to $M_{\Phi(\alpha)}$. 

Conversely, if $\widetilde f_n$ does not converge to a single orbit in $\xtilde_p$, let $\xi\in\xtilde_p$ be its limit point, which contains at least two non-trivial orbits. Because of Theorem \ref{theorem-MainTheorem}, the only possible norm limit
of $M_n$ is $M_\xi$. However, the interaction term between the two non-trivial pieces of mass which appears in the decomposition of $\Phi(f_n)$ described above, disallows convergence in norm to $M_\xi$.
\eproof

\noindent{\bf Proofs of Corollary \ref{corollary-limitpoints} and Corollary \ref{corollary-minimizers}.} These two results follow immediately from Lemma \ref{lemma-ContinuityPhi} combined with Theorem \ref{theorem-CompactificationMp}. \eproof

\section{Applications to amenability}\label{section-applications}

As mentioned in the introduction, our main motivation for studying the Schur multipliers appearing in Theorem \ref{theorem-MainTheorem} is that they naturally capture properties of the group. The main goal of this section is to state and prove two characterizations of amenability in terms of the Schur multipliers in $\bS_p = \bS_p(\Gamma)$ with $p\in[1,2)$ for a countable discrete group $\Gamma$.

\subsection{Amenability: Compact approximations and a variational characterization.}\label{sec-application1}

For our first main application, the idea is to approximate the non-compact space $\bS_2$, which captures amenability, by the compact spaces $\bX_p$. To that end, we introduce some additional notations. Define
\vspace{1mm}
$$
\bX \coloneqq \bX(\Gamma) \coloneqq \prod_{p\in[1,2)} \bX_p(\Gamma),
$$
which is compact with respect to the product topology. The image of $\xtilde$ in $\bX$ under the natural inclusion 
$\xi \mapsto (M_{\xi^{1/p}})_{p\in[1,2)}$ is a compact subset and will be denoted by $\bK \coloneqq \bK(\Gamma)$. For a finite subset $F$ of $\Gamma$ and a Schur mutliplier $M_k$, where $k$ is a kernel, we define $r_F(M_k)$ to be the multiplier {\it cut-off outside $\Tube(F)$}, i.e.\ we set
$$
r_F(M_k) \coloneqq M_{k \1_{\Tube(F)}}.
$$
Note that for Schur multipliers cut-off outside finite tubes, convergence in strong operator topology as well as convergence in norm coincide with pointwise convergence of the associated kernels on $\Tube(F)$.

\begin{theorem} \label{theorem: CompactApproximationId} Let $\Gamma$ be a countable discrete group and let $A$ denote either $C_\lambda^*(\Gamma)$ or $C_u^*(\Gamma)$. Then $\Gamma$ is amenable if and only if one of the following equivalent conditions hold:
\begin{enumerate}
\item The identity lies in the closure of 
$$
\bigcup_{p\in[1,2)} \bS_p(\Gamma)
$$
w.r.t.\ to the strong operator topology on $B(A)$.
\item There exist $\bx_n=(M_p^{(n)})_{p\in[1,2)} \in \bK(\Gamma)$, $n\in\N$, such that a diagonal sequence $(M_{p_k}^{(n_k)})_{k=1}^\infty$ converges to the identity w.r.t.\ to the strong operator topology on $B(A)$.
\item For each finite subset $F$ of $\Gamma$
\begin{equation} \label{equation: deltaFp}
\delta^{(F,p)} \coloneqq \min_{M \in \bX_p(\Gamma)} \Vert M_{\1_{\Tube(F)}} - r_F(M) \Vert 
\end{equation}
converges to zero as $p \in [1,2)$ converges to $2$. \\
\noindent In this case, for any sequence of finite sets $F_n \uparrow \Gamma$, there exist $p_n \in [1,2)$ such that the minimizers $M^{(F_n,p_n)}$ in (\ref{equation: deltaFp}) converge, as $n\to\infty$, to the identity w.r.t.\ to the strong operator topology on $B(A)$.
\end{enumerate}
\end{theorem}

\begin{remark} The minimum in (\ref{equation: deltaFp}) is guaranteed by the compactness proved in Theorem \ref{theorem-MainTheorem}. We point out that the operator norm may be replaced by any distance which metrizes the topology on the set of Schur multipliers cut-off outside $\Tube(F)$. An interpretation of this condition is the following: Given a distance and a fixed $p\in[1,2)$, there exists some $M\in\bS_p$ which minimizes the distance to the identity restricted to the set $F$ (recall that in $\bS_2$ there in general is no such minimizer). Letting $p_n \to 2$, the minimizers become optimal restricted to $F$ and letting $F_n\to\Gamma$ and $p_n \to 2$ simultaneously in a suitable way, the minimizers actually give us the desired approximation of the identity on all of $\Gamma$.
\end{remark}

\subsubsection{Background on amenability.} Let $\Gamma$ be a countable discrete group. An {\it invariant mean on } $\ell^\infty(\Gamma)$ is a state on $\ell^\infty(\Gamma)$ which is invariant under the left translation action, i.e.\ a positive linear functional $\mu$ on $\ell^\infty(\Gamma)$ with $\mu(1)=1$ and $\mu(L_s f)=\mu(f)$ for every $s\in\Gamma$, where $L_s f(x)=f(s^{-1}x)$ for $f\in\ell^\infty(\Gamma)$ and $s,x\in\Gamma$.
The group is called {\it amenable} if there exists an invariant mean on $\ell^\infty(\Gamma)$. For future reference and the convenience of the reader, we collect some well-known characterizations of amenability which will be used in this section.

\begin{prop} \label{prop: amenability} 
Let $\Gamma$ be a countable discrete group. Then the following are equivalent:
\begin{enumerate}
\item $\Gamma$ is amenable.
\item There exists a net of finitely supported positive definite functions converging pointwise to $1$.
\item There exists a F{\o}lner sequence, that is a sequence $(F_n)_{n=1}^\infty$ of non-empty finite subsets of $\Gamma$ such that
$$
\lim_{n\to\infty} \frac{|F_n \Delta sF_n|}{|F_n|}=0 \quad \text{for every } s\in \Gamma.
$$

\item There exists a right-F{\o}lner sequence, that is a sequence $(F_n)_{n=1}^\infty$ of non-empty finite subsets of $\Gamma$ such that
$$
\lim_{n\to\infty} \frac{|F_n \Delta F_n s|}{|F_n|}=0,
$$
or equivalently
$$
\lim_{n\to\infty} \frac{|F_n \cap F_ns|}{|F_n|}=1 \quad \text{for every } s\in \Gamma.
$$
\item There exists a sequence of almost right-invariant unit vectors in $\ell^2(\Gamma)$, that is a sequence $(f_n)_{n=1}^\infty$ of unit vectors in $\ell^2(\Gamma)$ such that
$$
\lim_{n\to \infty} \Vert f_n - f_n * \delta_s \Vert_2^2 = 0 \quad \text{for every } s\in \Gamma.
$$
\end{enumerate}
\end{prop}

\proof 
These are standard facts, for proofs of which we refer the reader to Section 2.6 and in particular Theorem 2.6.8 in \cite{BrownOzawa}. Changing between left- and right-invariant statements is straightforward by taking inverses in the group. \eproof

\subsubsection{Proof of Theorem \ref{theorem: CompactApproximationId}.} 
Let $\Gamma$ be a countable discrete group. We start this subsection with three preliminary observations and then proceed to prove Theorem \ref{theorem: CompactApproximationId}. We start by showing that the property of almost right-invariance is characterized by the image of the sequence under $\Phi$, i.e.\ by the corresponding positive definite functions.

\begin{lemma} \label{lemma: ConvPhiTo1} Let $(f_n)_{n=1}^\infty$ be in $\Mcal_2$. Then the sequence $(f_n)_n$ is almost right-invariant if and only if the sequence $(\Phi(f_n))_n$ converges pointwise to $1$.
\end{lemma}

\proof For every $s \in \Gamma$ we may rewrite
\begin{equation} \label{equ: EqualityPhi}
\Vert f_n - f_n * \delta_s \Vert_2^2 = \langle (1-\rho(s))f_n,(1-\rho(s))f_n \rangle = 2\bigl( 1-\langle \rho(s)f_n, f_n \rangle \bigr) = 2\bigl(1-\Phi(f_n)(s)\bigr),
\end{equation}
which implies the claim. \eproof

The next result is related to the well-known characterization of amenability via the existence of an approximate identity for the Fourier algebra. We include a proof for the sake of completeness.

\begin{lemma}[Approximate Identity Criterion] \label{lemma: amenabilitycriterionM2}  Let $A$ denote either $C_\lambda^*(\Gamma)$ or $C_u^*(\Gamma)$. Then $\Gamma$ is amenable if and only if the identity lies in the closure of $\bS_2(\Gamma)$ with respect to the strong operator topology on $B(A)$.
\end{lemma}

\proof Suppose that $\Gamma$ is amenable and let $(F_n)_{n=1}^\infty$ be a right F{\o}lner sequence. Then $f_n \coloneqq |F_n|^{-1/2} \1_{F_n}$ defines a sequence in $\Mcal_2$ such that for every $s \in \Gamma$
$$
\lim_{n \to \infty} \Phi(f_n)(s) = \lim_{n \to \infty} \frac{1}{|F_n|} \sum_{x \in \Gamma} \1_{F_n}(x) \1_{F_n}(xs) = \lim_{n \to \infty} \frac{|F_n \cap F_n s^{-1}|}{|F_n|} = 1,
$$
compare Proposition \ref{prop: amenability} (iv). Now Lemma \ref{lemma-strongconvergence} implies 
$$
\lim_{n \to \infty} M_{\Phi(f_n)} = M_{\mathbf 1}
$$
in the strong operator topology on $B(A)$. To prove the converse, we choose a sequence in $\bS_2$ converging to the identity in the strong operator topology on $B(A)$. By Lemma \ref{lemma-strongconvergence} we obtain a sequence $(f_n)_{n=1}^\infty \subset \Mcal_2$ such that
$$
\lim_{n \to \infty} \Phi(f_n) (x) = 1 \quad \text{for every } x \in \Gamma.
$$
Applying Lemma \ref{lemma: ConvPhiTo1} to this sequence we conclude that $(f_n)_n$ is a sequence of almost right-invariant unit vectors in $\ell^2(\Gamma)$ and hence $\Gamma$ is amenable by Proposition \ref{prop: amenability} (v). \eproof

\begin{lemma} \label{lemma: l2ContinuityPhi}
Let $(f_n)_{n=1}^\infty$ be in $\ell^2(\Gamma)$ such that $f_n \to f \in \ell^2(\Gamma)$. Then $\Phi(f_n)\to\Phi(f)$ pointwise.
\end{lemma}

\proof For every $s\in\Gamma$ the Cauchy Schwarz inequality implies that
\begin{flalign*}
| \Phi(f)(s) - \Phi(f_n)(s) | & = \bigl| \langle f,\rho(s)f \rangle - \langle f_n,\rho(s)f_n \rangle \bigr| \\
& \leq \bigl| \langle f,\rho(s)f \rangle - \langle f_n,\rho(s)f \rangle \bigr| + \bigl| \langle f_n,\rho(s)f \rangle - \langle f_n,\rho(s)f_n \rangle \bigr| \\
& \leq \Vert f-f_n \Vert_2 \Vert \rho(s)f \Vert_2 + \Vert f_n \Vert_2 \Vert \rho(s)(f-f_n) \Vert_2 \\
& \leq \Bigl( \Vert f \Vert_2 + \sup_{n\in\N} \Vert f_n \Vert_2 \Bigr) \Vert f-f_n \Vert_2,
\end{flalign*}
which converges to zero as $n\to\infty$. \eproof

\noindent {{\bf Proof of Theorem \ref{theorem: CompactApproximationId}}. Throughout this proof, {\it convergence in strong operator topology} will refer to convergence in the strong operator topology on $B(A)$. We first show that (i) holds if and only if $\Gamma$ is amenable.  Assume that $\Gamma$ is amenable, then, by Lemma \ref{lemma: amenabilitycriterionM2}, there exists a sequence in $\bS_2$ converging to the identity in strong operator topology. In order to prove that (i) holds, it therefore suffices to show that every $M \in \bS_2$ can be written as a strong operator topology limit of a sequence in 
$$
\bigcup_{p\in[1,2)} \bS_p.
$$
To see this, fix some $M_{\Phi(f)}$ with $f\in\Mcal_2$. Choose $p_n \in [1,2)$ such that $p_n \to 2$ and consider the sequence $(g_n)_{n=1}^\infty$ defined by $g_n \coloneqq f^{2/p_n} \in \Mcal_{p_n}$. Then $g_n \to f$ in $\ell^2(\Gamma)$ because of pointwise convergence and convergence of the norms. Hence Lemma \ref{lemma: l2ContinuityPhi} implies that $\Phi(g_n) \to \Phi(f)$ pointwise. It follows from Lemma \ref{lemma-strongconvergence} that $M_{\Phi(g_n)} \in \bS_{p_n}$ converges to $M_{\Phi(f)}$ as in strong operator topology.

Conversely, suppose that there exist $p_n \in [1,2)$ and operators $M_{p_n}^{(n)} \in \bS_{p_n}$ which converge to the identity in strong operator topology. Write $M_{p_n}^{(n)} = M_{\Phi(f_n)}$ for suitable $f_n\in\Mcal_{p_n}$, then Lemma \ref{lemma-strongconvergence} shows that the sequence $(\Phi(f_n))_n$ converges pointwise to $1$. Now note that each $f_n$ can be approximated in $\ell^2$-norm by finitely supported functions. By Lemma \ref{lemma: l2ContinuityPhi} and Lemma \ref{lemma: PhifPosDef}, after taking a suitable diagonal sequence, we find a sequence of finitely supported positive definite functions converging to $1$. This proves amenability of $\Gamma$ by Proposition \ref{prop: amenability} (ii).

Clearly (i) and (ii) are equivalent. Finally, we verify that the conditions (i) and (iii) are equivalent.

(i) $\Rightarrow$ (iii): We have already seen that (i) implies amenability of $\Gamma$. Assuming that (i) holds, we may therefore choose a right F{\o}lner sequence $(F_n)_{n=1}^\infty$. Choose a sequence $(p_n)_{n=1}^\infty$ in $[1,2)$ with $p_n\to2$ such that $|F_n|^{1-2/p_n} \to 1$ as $n\to\infty$ and define $f_n \coloneqq |F_n|^{-1/p_n} \1_{F_n} \in \Mcal_{p_n}$. The same computation as in Lemma \ref{lemma: ConvPhiTo1} shows that for every $s\in\Gamma$
\begin{equation} \label{equation: CompApproxEqu1}
\begin{aligned}
\Phi(f_n)(s) = \Vert f_n \Vert_2^2 - \frac{1}{2} \, \Vert f_n - \rho(s) f_n \Vert_2^2 &= |F_n|^{1-2/p_n} - \frac{|F_n \Delta F_n s^{-1}|}{|F_n|^{2/p_n}} 
\\
&\geq |F_n|^{1-2/p_n} - \frac{|F_n \Delta F_n s^{-1}|}{|F_n|},
\end{aligned}
\end{equation}
where the last inequality is due to the fact that $p_n \leq 2$. By the F{\o}lner property and our assumption on $p_n$, the right-hand side of Inequality (\ref{equation: CompApproxEqu1}) converges to $1$ as $n\to\infty$. This proves that the sequence $(\Phi(f_n))_n$ converges pointwise to $1$, hence uniformly on the finite set $F$ and thus by Lemma \ref{lemma: FinPropBound}
$$
\lim_{n\to\infty} \delta^{(F,p_n)} \leq \lim_{n\to\infty}  \Vert M_{\1_{\Tube(F)}} - r_F(M_{\Phi(f_n)}) \Vert \leq \lim_{n\to\infty} |F| \, \Vert (1-\Phi(f_n))_{|F} \Vert_\infty = 0,
$$
which proves the claim.

(iii) $\Rightarrow$ (i): Let $F_n \subset \Gamma$ finite such that $F_n \uparrow \Gamma$ and let $\eps_n>0$ with $\eps_n \downarrow 0$. Condition (\ref{equation: deltaFp}) and the fact that $\bS_p \subset \bX_p$ is dense imply that we find for every $n$ some $f_{p_n}\in\Mcal_{p_n}$ such that
$$
\Vert M_{\1_{\Tube(F_n)}}-r_F(M_{\Phi(f_{p_n})}) \Vert \leq \eps_n.
$$
In particular, we must have $|1-\Phi(f_{p_n})(s)| \leq \eps_n$ for every $s\in F_n$. By our assumptions on $F_n$ and $\eps_n$, $\Phi(f_{p_n})$ converges pointwise to $1$. Lemma \ref{lemma-strongconvergence} yields that $M_n \coloneqq M_{\Phi(f_{p_n})}$ defines a sequence in 
$$
\bigcup_{p\in[1,2)} \bS_p
$$
converging to the identity in strong operator topology. Note that the same argument clearly shows convergence of the minimizers as well. The proof of Theorem \ref{theorem: CompactApproximationId} is thus complete. \eproof

\subsection{Schur multipliers and group-invariant percolation.}\label{sec-percolation}

We turn to our second characterization of amenability, for which we first set the background of {\it group-invariant percolation}.
Let $\Gamma$ be a finitely generating group and let $S$ be a finite and symmetric generating set. Then the {\it (right) Cayley graph} of $\Gamma$ with respect to $S$, which will be denoted $X(\Gamma,S)$, is the undirected graph with vertex set $\Gamma$ and an edge between $x,y\in\Gamma$ if and only if $x^{-1}y\in S$. Clearly this graph depends on the generating set, however, it is easy to see that two Cayley graphs of $\Gamma$ with respect to different finite and symmetric generating sets are quasi-isometric. We note that $\Gamma$ has a canonical action on $X(\Gamma,S)$ via left-multiplication. 

Given any countable, locally finite graph $X=(V,E)$, we endow the set $\{0,1\}^V$ with the Borel-$\sigma$-field generated by the cylinder events. An element $\omega \in \{0,1\}^V$ is called a {\it configuration} and is identified with the subgraph of $X$ which has vertex set 
$$
V(\omega) = \{ x\in V \colon \omega_x = 1\}
$$
 and edge set $E(\omega)$ consisting of all edges in $X$ which connect vertices in $V(\omega)$. For a vertex $x$, the {\it cluster} of $x$ is the connected component of $x$ in the configuration and is denoted $K(x)$. A probability measure on the set of configurations is called a {\it site percolation} on $X$. Given a subgroup $G$ of the automorphism group ${\rm Aut}(X)$ of $X$, a site percolation is called {\it $G$-invariant} if it is invariant under the induced action of $G$. Given the above background, we can state our next main result: 

\begin{theorem} \label{theorem: PercolationConstruction} Let $\Gamma$ be a finitely generated group, let $X$ be one of its Cayley graphs w.r.t.\ some finite and symmetric generating set and let $(p_n)_{n=1}^\infty \subset [1,2)$ such that $p_n \to 2$. Then $\Gamma$ is amenable if and only if there exists a sequence $(\P_n)_{n=1}^\infty$ of $\Gamma$-invariant site percolations on $X$ such that the following two conditions hold:$^{\mathrm d}$\footnote{$^{\mathrm d}$For a probability measure $\P$ (resp. $\P_n$), we denote by $\E=\E^\P$ (resp. $\E_n=\E^{\P_n}$) the corresponding expectation.}
\begin{itemize}
\item Each $\P_n$ has no infinite components, and 
\item For every $s\in\Gamma$
$$
\lim_{n\to\infty} \E_n\left[ \Phi \left( |K(o)|^{-1/p_n} \1_{K(o)} \right)(s) \right] = 1,
$$
where $K(o)$ denotes the random cluster of the identity in $X$ and 
$$
|K(o)|^{-1/p_n} \1_{K(o)} \coloneqq 0\quad\mbox{ if }K(o)=\emptyset.
$$
\end{itemize}
Moreover, in this case 
$$
\varphi_n(s)\coloneqq\E_n\left[ \Phi \left( |K(o)|^{-1/p_n} \1_{K(o)} \right)(s) \right]
$$
 is positive definite and the Schur multipliers $M_{\varphi_n}$ on $B(C_\lambda^*(\Gamma))$, respectively on $B(C_u^*(\Gamma))$, converge to the identity in the strong operator topology as $n\to\infty$.
\end{theorem}

\begin{remark} 
It is not difficult to show that for every $\Gamma$-invariant site percolation $\P$ on $X$ the function
$$
\varphi(s) = \E\left[ \Phi \left( |K(o)|^{-1/p_n} \1_{K(o)} \right)(s) \right]
$$
is still positive definite. Hence the corresponding sequence $M_{\varphi_n}$ arising from Theorem \ref{theorem: PercolationConstruction} provides an interesting probabilistic example of the type of Schur multipliers we have been considering so far.
\end{remark}

\begin{remark} 
The assumption that $\Gamma$ be finitely generated is needed in order to work with the Cayley graph of $\Gamma$. Many interesting countable discrete groups are finitely generated, so this should not be considered a strong restriction.
\end{remark}

 \subsubsection{Proof of Theorem \ref{theorem: PercolationConstruction}.} 

We start by recalling the following celebrated result of Benjamini, Lyons, Peres and Schramm, which provides a general characterization of amenability for automorphism groups of graphs.

\begin{theorem}[{\cite[Theorem 5.1]{BLPS99}}] \label{theorem: BLPSCriterion} Let $X$ be a countable, locally finite graph and $G$ a closed subgroup of ${\rm Aut}(X)$. Then $G$ is amenable if and only if for every $\alpha<1$, there exists a $G$-invariant site percolation $\P$ on $X$ with no infinite components and with $\P[x\in\omega]>\alpha$ for all vertices $x$.
\end{theorem}

We will in the following consider the special case of Theorem \ref{theorem: BLPSCriterion} when $G$ is a finitely generated group acting on one of its Cayley graphs via left-multiplication. We will now prove that the invariant percolations appearing in the above theorem are also of the type appearing in Theorem \ref{theorem: PercolationConstruction}.

\noindent {\bf Proof of Theorem \ref{theorem: PercolationConstruction}.} 
First assume that there exists a sequence $(\P_n)_{n=1}^\infty$ of $\Gamma$-invariant site percolations on $X$ such that each $\P_n$ has no infinite components and for every $s\in\Gamma$
$$
\lim_{n\to\infty} \E_n\left[ \Phi \left( |K(o)|^{-1/p_n} \1_{K(o)} \right)(s) \right] = 1.
$$
Set $f_{p_n} = |K(o)|^{-1/p_n} \1_{K(o)}$ if $K(o)\neq\emptyset$ and $f_{p_n} = 0$ otherwise. Then since $p_n \leq 2$
$$
\Phi(f_{p_n})(o) = \Vert f_{p_n} \Vert_2^2 = |K(o)|^{1-2/p_n} \1_{\{K(o) \neq \emptyset\}} \leq \1_{\{o\in\omega\}} 
$$
and thus
$$
\E_n[\Phi(f_{p_n})(o)] \leq \P_n[o\in\omega].
$$
Therefore 
$$
\liminf_{n\to\infty} \P_n[o\in\omega] \geq \lim_{n\to\infty} \E_n\left[ \Phi \left( f_{p_n} \right)(o) \right] = 1
$$
and amenability of $\Gamma$ follows from Theorem \ref{theorem: BLPSCriterion}.

\noindent We prove the converse in two steps. 

\noindent{\bf Step 1:} We first show that if $\Gamma$ is amenable, then there exists a sequence $(\P_n)_{n=1}^\infty$ of $\Gamma$-invariant site percolations on $X$ such that each $\P_n$ has no infinite components and for every $s\in\Gamma$
$$
 \lim_{n\to\infty} \E_n\left[ \Phi \left( f_2 \right)(s) \right] = 1,
$$
where of course $f_2 = |K(o)|^{-1/2} \1_{K(o)}$ if $K(o)\neq\emptyset$ and $f_2 = 0$ otherwise. 

To see this, use Theorem \ref{theorem: BLPSCriterion} to find a sequence $(\P_n)_{n=1}^\infty$ of $\Gamma$-invariant site percolations on $X$ such that each $\P_n$ has no infinite components and $\P_n[o\in\omega] \to 1$ as $n\to\infty$. Note that if $K(o) \neq \emptyset$, we have the formula
$$
\Phi(f_2)(s)=1-\frac{1}{2}\, \Vert f_2 - \rho(s)f_2 \Vert_2^2,
$$
compare Equation (\ref{equ: EqualityPhi}). Hence
$$
\E_n[\Phi(f_2)(s)] = \P_n[o \in \omega] - \frac{1}{2} \, \E_n \Bigl[ \Vert f_2 - \rho(s)f_2 \Vert_2^2 \Bigr]
$$
and it suffices to prove that the second term vanishes as $n\to\infty$. To see this, we rewrite 
$$
\E_n \Bigl[ \Vert f_2 - \rho(s)f_2 \Vert_2^2 \Bigr] = \E \biggl[ \frac{|K(o)\setminus K(o)s^{-1}|}{|K(o)|} \biggr] + \E \biggl[ \frac{|K(o)s^{-1}\setminus K(o)|}{|K(o)|} \biggr].
$$
Both summands can be estimated similarly using the mass-transport principle  \cite[Eq.\ (2.1) or p.\ 15]{BLPS99}, which states that for any diagonally $\Gamma$-invariant function $f\colon V\times V \to [0,\infty]$, we have that $\sum_{y\in V} f(y,o)=\sum_{x\in V} f(o,x)$.
We do this for the first: Define the following invariant mass-transport: Put a unit mass at every vertex which lies in a finite component but not its $s^{-1}$-right-translate, then redistribute uniformly over this component; or more precisely
$$
m(x,y,\omega) = \1\bigl\{ x \in \omega, xs \notin K(x),y \in K(x) \bigl\} \frac{1}{|K(x)|}.
$$
The mass-transport principle  applied to $f(x,y):=\E_n[m(x,y,\omega)]$ implies that
$$
\E_n \biggl[ \frac{|K(o)\setminus K(o)s^{-1}|}{|K(o)|} \biggr] = \P_n[o \in \omega, s \notin K(o)].
$$
Choose a finite path $\rho_s = \{x_1,\ldots,x_N\}$ connecting $o$ to $s$. Since $\P_n$ is invariant and $\P_n[o\in\omega]\to1$, it follows that
$$
\P_n[\rho_s {\rm \ is \ contained \ in \ } K(o)] = \P_n[o\in\omega, x_i \in \omega {\rm \ for \ all \ } i=1,\ldots,N] \to 1 \quad \text{as } n \to \infty
$$
because we have a finite intersection of events with equal probability close to $1$. Therefore
$$
\P_n[o \in \omega, s \notin K(o)]  \leq \P_n[ \rho_s {\rm \ is \ not \ contained \ in \ } K(o)] \to 0 \quad \text{as } n \to \infty.
$$
This concludes the proof of Step 1.

\noindent{\bf Step 2:} Let $\P_n$ be a sequence as in Step 1. For $t\geq1$, the function $p \mapsto t^{-1/p}$ is non-decreasing, hence $f_{p_n}$ converges monotonically to $f_2$ and thus $\Phi(f_{p_n})$ converges monotonically to $\Phi(f_2)$. Therefore the monotone convergence theorem implies that for every $s\in\Gamma$
\vspace{1mm}
\begin{equation} \label{eq-monConv}
\E_n[\Phi(f_p)(s)] \uparrow \E_n[\Phi(f_2)(s)] \quad \text{as } p \uparrow 2.
\end{equation}
Let $F\subset\Gamma$ be a finite subset and $\eps>0$. By Step 1, there exists $n^{(F,\eps)}\in\N$ such that 
$$
\E_{n^{(F,\eps)}} [\Phi(f_2)(s)] \geq 1- \eps/2 \quad \mbox{for all} \, \, s \in F.
$$
By (\ref{eq-monConv}), there exists $p^{(F,\eps)} \in [1,2)$ such that
$$
\E_{n^{(F,\eps)}} [\Phi(f_p)(s)] \geq 1- \eps \quad \mbox{for all} \, \, s \in F \, \, \mbox{and} \, \,  p\geq p^{(F,\eps)}.
$$
Now choose an increasing sequence of finite subsets $F_i \uparrow \Gamma$ and let $\eps_i \in(0,1)$ such that $\eps_i \downarrow 0$. For each $i\in\N$, let $n_i^* \coloneqq n^{(F_i,\eps_i)}$ and $p_i^* \coloneqq p^{(F_i,\eps_i)}$ be as defined above. Without loss of generality we may assume that $(p_i^*)_{i=1}^\infty \subset [1,2)$ forms an increasing sequence. Moreover, we necessarily have $p_i^*\uparrow2$. Since $p_n\to2$, it follows that for every $i\in\N$ there exists $n_i\in\N$ such that $p_n \geq p_i^*$ for all $n\geq n_i$. Without loss of generality $(n_i)_{i=1}^\infty$ forms an increasing sequence. For each $n\in\N$, we then have $n \in \{n_i,n_{i+1}\}$ for some $i\in\N$ and we define $\mathbf P_n=\P_{n_i^*}$. Note that each $\mathbf P_n$ is a $\Gamma$-invariant site percolation with no infinite components because it is chosen from the sequence $(\P_n)_n$ from Step 1. Moreover, if $n\in\{n_i,n_{i+1}\}$ we have that $p_n\geq p_i^*$ and thus (denoting expectation w.r.t. $\mathbf P_n$ by $\mathbf E_n$),
$$
\mathbf E_n [\Phi(f_{p_n})(s)] = \E_{n_i^*} [\Phi(f_{p_n})(s)] \geq \E_{n_i^*} [\Phi(f_{p_i^*})(s)] \quad \mbox{for all} \, \, s \in \Gamma.
$$
Finally, for every $s\in\Gamma$ we have that $s\in F_i$ for all large enough $i$ and therefore
$$
\liminf_{n\to\infty} \mathbf E_n [ \Phi(f_{p_n})(s)] \geq \liminf_{i\to\infty} \E_{n_i^*} [ \Phi(f_{p_i^*})(s)] \geq \lim_{n \to \infty} 1- \eps_i = 1.
$$
Thus the sequence $(\mathbf P_n)_{n=1}^\infty$ has the desired properties. 

Finally, note that 
$$
s\mapsto   \Phi \left( |K(o)|^{-1/p_n} \1_{K(o)} \right)(s) 
$$
is positive definite. Taking expectation preserves positive-definiteness, thus showing that 
$
\varphi(s) = \E\left[ \Phi \left( |K(o)|^{-1/p_n} \1_{K(o)} \right)(s) \right]
$
is also positive definite, proving Theorem \ref{theorem: PercolationConstruction}. 
\eproof

\noindent{\bf Acknowledgement:} It is a pleasure to thank Mahan Mj. (Mumbai) who pointed out to the first author the theory of Roe algebras and suggested to investigate its links to \cite{MV16} that inspired this work. 
The authors would like to thank an anonymous referee for a careful reading and suggesting a number of improvements. 
The research of both authors is funded by the Deutsche Forschungsgemeinschaft (DFG) under Germany's Excellence Strategy EXC 2044-390685587, Mathematics M\"unster: Dynamics-Geometry-Structure.


\begin{thebibliography}{WWW98}

\bibitem{BLPS99} 
{\sc Benjamini, I., Lyons, R., Peres, Y., Schramm, O.} (1999). Group-invariant percolation on graphs. {\em Geom. and Funct. Anal.} \textbf{9} 29-66.

\bibitem{B77}
{\sc Bennet, G.} (1977). Schur multipliers. \textit{Duke Math. J.} \textbf{44} 603-639.

\bibitem{Billingsley}
{\sc Billingsley, P.} Convergence of probability measures. In: Wiley series in probability and statistics. Second Edition, John Wiley \& Sons, Inc., New York, 1999.

\bibitem{BrownOzawa}
{\sc Brown, N. P.} and {\sc Ozawa, N.} $C^*$-Algebras and finite dimensional approximations. In: Graduate Studies in Mathematics, Vol. \textbf{88}, AMS, Providence, Rhode Island, 2008.

\bibitem{CW04}
{\sc Chen, X.} and {\sc Wang, Q.} (2004). Ideal structure of uniform Roe algebras of coarse spaces. \textit{J. Funct. Anal.} \textbf{216} 191-211.

\bibitem{DD05}
{\sc Davidson, D.} and {\sc Donsig, A. P.} (2005). Norms of Schur multipliers. \textit{Illinois J. Math.} \textbf{51} 743-766.

\bibitem{DK18}
{\sc Dym, H.} and {\sc Katsnelson, V.} Contributions of Issai Schur to analysis. In: Studies in memory of Issai Schur (Chevaleret/Rehovot, 2000), Progr. Math. \textbf{210}, Birkhäuser Boston, Boston, MA,  2003.

\bibitem{E64}
{\sc Eymard, P.} (1964). L'algèbre de Fourier d'un groupe localement compact. \textit{Bull.\ S.\ M.\ F.} \textbf{92} 181-236.

\bibitem{Grothendieck}
{\sc Grothendieck, A.} (1953). Résumé de la théorié métrique des produits tensoriels topologiques. \textit{Bull. Soc. Mat. $S\tilde{a}o$-Paulo} \textbf{8} 1-79.

\bibitem{GK02}
{\sc Guentner, E.} and {\sc Kaminker, J.} (2002). Exactness and the Novikov conjecture. \textit{Topology} \textbf{41}(2) 411-418.

\bibitem{Haagerup78}
{\sc Haagerup, U.} (1979). An example of a non nuclear $C^*$-Algebra, which has the metric approximation property. \textit{Invent. Math.} \textbf{50} 279-293.

\bibitem{HSS10}
{\sc Haagerup, U.}, {\sc Steenstrup, T.}, {\sc Szwarc, R.} (2010). Schur multipliers and spherical functions on homogeneous trees. \textit{Internat. J. Math.} \textbf{21}(10) 1337-1382.

\bibitem{H97}
{\sc Häggström, O.} (1997). Infinite clusters in dependent automorphism invariant percolation on trees. \textit{Ann. Prob.} \textbf{25} 1423-1436.

\bibitem{MV16}
{\sc Mukherjee, C.} and {\sc Varadhan, S. R. S.} (2016). Brownian occupation measures, compactness and large deviations. \textit{Ann. Prob.} \textbf{44} 3934-3964.

\bibitem{MR23}
{\sc Mukherjee, C.} and {\sc Recke, K.}, (2023). Haagerup property and group-invariant percolation. {\em arxive-preprint:} {\url{https://arxiv.org/abs/2303.17429}}.

\bibitem{O00}
{\sc Ozawa, N.} (2000). Amenable actions and exactness for discrete groups. \textit{C. R. Acad. Sci. Paris. Sér. I Math.} \textbf{330}(8) 691-695.

\bibitem{PRV62}
{\sc Parthasarathy, K. R.}, {\sc Rao, R.} and {\sc Varadhan, S. R. S.}(1962). On the Category of Indecomposable Distributions on Topological Groups. \textit{Trans. Amer. Math. Soc.} {\bf 102}, 200-217

\bibitem{Pisier1996}
{\sc Pisier, G.} Similarity problems and completely bounded maps. In: Lecture notes in mathematics; 1618, Springer-Berlin-Heidelberg, 1996.

\bibitem{Roe}
{\sc Roe, J.} Lectures on coarse geometry. In: University Lecture Series, Vol. \textbf{31}, AMS, Providence, Rhode Island, 2003.

\bibitem{RW14}
{\sc Roe, J.} and {\sc Willett, R.} (2014). Ghostbusting and property A. \textit{J. Funct. Anal.} \textbf{266}(3) 1674-1684.


\bibitem{WZ18}
{\sc Winter, W.} and {\sc Zacharias, J.} (2018). The nuclear dimension of $C^*$-algebras. \textit{Adv. Math.} \textbf{224} 461-498.

\bibitem{Y00}
{\sc Yu, G.} (2000). The coarse Baum-Connes conjecture for spaces which admit a uniform embedding into Hilbert space. \textit{Invent. Math.} \textbf{139} 201-240.

\end{thebibliography}
\end{document}